\documentclass[10pt]{amsart}

\usepackage{epsfig}
\usepackage{amsthm}
\usepackage{amssymb}
\usepackage{amsmath}
\usepackage{epic}
\usepackage{eepic}
\usepackage{graphicx}
\usepackage{psfrag}

\newcommand     {\comment}[1]   {}
\newcommand{\mute}[2] {}
\newcommand     {\printname}[1] {}
\newtheorem {Theorem}   {Theorem}
\numberwithin{Theorem}{section}
\newtheorem {Lemma}[equation]    {Lemma}
\newtheorem {Claim}[equation]    {Claim}
\newtheorem {Proposition}[equation]{Proposition}
\theoremstyle{definition}
\newtheorem{Definition}[equation]{Definition}
\theoremstyle{remark}
\newtheorem{Remark}[equation]{Remark}
\newtheorem{Example}[equation]{Example}

\def    \bfC    {{\mathbf C}}
\def    \bfS    {{\mathbf S}}
\def    \bfQ    {{\mathbf Q}}

\def    \bfN    {{\mathbf N}}
\def    \bfone  {{\mathbf 1}}

\def    \C      {\mathbb{C}}
\def    \R      {\mathbb{R}}

\def    \Z      {\mathbb{Z}}

\def    \PR     {\mathcal{R}}
\def    \PB     {\mathcal{B}}

\def    \bvert     {\substack{\vert \\[-.5mm] \vert \\[-.5mm] \vert}}

\def   \Cinf   {C^\infty}

\def    \codim    {{\operatorname{codim}}}

\newcommand{\dsum}{\displaystyle\sum}  
\newcommand{\dprod}{\displaystyle\prod}

\begin{document}

\bibliographystyle{plain}

\title[New polytope decompositions and Euler-Maclaurin formulas]
{New polytope decompositions and Euler-Maclaurin formulas for
simple integral polytopes}

\author[J.\ Agapito]{Jos\'e Agapito}
\author[L.\ Godinho]{Leonor Godinho}
\address{DEPARTAMENTO DE MATEM\'ATICA, INSTITUTO SUPERIOR T\'ECNICO, AV. RO\-VIS\-CO PAIS,
1049-001 LISBON, PORTUGAL, FAX: (351) 21 841 7035}
\email{agapito@math.ist.utl.pt}
\email{lgodin@math.ist.utl.pt}
\thanks{The first author was  partially supported by FCT (Portugal) through program POCTI/FEDER and grant POCTI/SFRH/BPD/20002/2004; the second author was partially supported by FCT through program POCTI/FEDER and  grant POCTI/MAT/57888/2004, and by Funda\c{c}\~{a}o Calouste Gulbenkian.}

\begin{abstract}
We use a version of localization in equivariant cohomology for the  norm-square  of the 
moment map, described by Paradan, to give several weighted decompositions  for simple polytopes.
As an application, we study Euler-Maclaurin formulas.
\end{abstract}

\maketitle


\section{Introduction}
\label{se:intro} The interplay between symplectic geometry and
combinatorics through the study of moment maps and the use of
equivariant cohomology and the geometry of toric varieties is a well
known and fertile theme in mathematics. See for example the work of Brion-Vergne \cite{BrV}, Cappell-Shaneson \cite{CS1}, \cite{CS2}, Ginzburg-Guillemin-Karshon \cite{GGK}, Guillemin \cite{Gu1}, \cite{Gu2},  Morelli \cite{M}, and Pommershein and Thomas \cite{PT}, \cite{Po}.

In this paper we use a version of localization in equivariant cohomology for the  norm-square  of the 
moment map, due to Paradan (see \cite{P}), to  motivate several new weighted polytope 
decomposition formulas. Indeed, in Section~\ref{sec:3}, applying this localization principle to the particular case of toric manifolds, we obtain weighted polytope decompositions for Delzant polytopes \cite{De} (see Example~\ref{ex:dec}) and then, by a purely combinatorial argument, we show, in Section~\ref{se:decompositions}, that they are in fact  valid   for \emph{any} simple\footnote{A polytope in $\R^d$ is called simple if each  vertex is the intersection of exactly $d$ facets
(i.e. codimension-$1$ faces).}  polytope not necessarily the moment map image of a toric manifold  (cf. Theorem~\ref{thm:1}). Moreover, still in Section~\ref{se:decompositions}, we use these decompositions to obtain new ones (Theorem~\ref{th:WBG+}) that generalize both the  Lawrence-Varchenko decomposition (see \cite{V} and \cite{L}) and the  Brianchon-Gram formula  (see \cite{Br}, \cite{G}, \cite{So} and \cite{Gr}) (cf. Remark~\ref{rmk:general}).

The well-known classical polytope decomposition formula of  Brianchon-Gram  expresses  the characteristic function $\bfone_P$ of a convex polytope $P$ as the alternating sum of the characteristic functions of all tangent cones to the faces of $P$. By flipping the edge vectors emanating from each vertex of $P$ in a systematic way using a polarizing vector, we obtain the Lawrence-Varchenko decomposition (also known as polar decomposition) which expresses the characteristic function of a convex simple polytope (\emph{only}) in terms of the characteristic functions of polarized cones supported at the vertices. Karshon, Sternberg and Weitsman \cite{KSW} and Agapito \cite{A03} gave weighted versions of this decomposition by assigning  weights to the faces of the polytope and of the cones in a consistent way. Our polytope decompositions combine the above two formulas. Like Brianchon-Gram they express   $\bfone_P$ in terms of characteristic functions of cones with apex the different faces of the polytope. However, these cones may no longer be the ordinary tangent cones to the polytope. Indeed, in our formulas, to each face is assigned a \emph{different} polarizing vector and we flip the edges of the tangent cones accordingly. In the first decomposition formula (Theorem~\ref{thm:1}) these polarizing vectors are obtained by choosing a suitable starting point $\varepsilon$ (the same for all vectors) and then taking as end points its orthogonal projections $\beta(\varepsilon,F)$ onto the  faces $F$ of the polytope, whenever these projections are nonempty (cf. Figure~\ref{fi:acosta1}). In the second decomposition formula  (Theorem~\ref{th:WBG+}) we take the vectors $\varepsilon - \beta(\varepsilon,F)$ as polarizing vectors (instead of $\beta(\varepsilon,F) -\varepsilon$ as above). This second formula  generalizes both the  Lawrence-Varchenko and the  Brianchon-Gram relations in the sense that, choosing $\varepsilon$ in suitable regions of $\R^d$, we obtain these polytope decompositions.      

As an application, in Section~\ref{se:EML}, we use the new  decompositions to give new Euler-Maclaurin formulas with remainder similar to those of Karshon-Sternberg-Weitsman \cite{KSW} and Agapito-Weitsman \cite{AW}. The classical Euler-Maclaurin formula computes the sum of the values of a function $f$  over the integer points of an interval in terms of the integral of $f$ over variations of that interval. This formula was generalized by Khovanskii and Pukhlikov (see \cite{KP1} and  \cite{KP2})  to a formula for the sum of the values of an exponential or polynomial function on the lattice points of a regular integral polytope, by Cappell and Shaneson \cite{CS1}, \cite{CS2}, \cite{S}, Guillemin \cite{Gu2} and Brion and Vergne \cite{BrV},  to \emph{simple} integral polytopes, and by Berline, Brion, Szenes and Vergne \cite{BV2}, \cite{BrV2} \cite{SV}, to any rational polytope. Note that all these formulas are exact and valid for sums of values of exponential or polynomial functions. Moreover, the formula in \cite{BV2} has the additional feature that it is \emph{local}, in the sense that it is given as a sum of integrals over the faces of the polytope of maps $D(F) \cdot p$, for operators $D(F)$ depending \emph{only} on a neighborhood of a generic point of the face $F$. In \cite{KSW}, Karshon, Sternberg and Weitsman prove a formula with remainder for the sum of an arbitrary smooth function $f$ of compact support, on  the integer points of a simple polytope. There, the remainder is given as a sum over the vertices of the polytope, of integrals over cones with those vertices,  of bounded periodic functions times several partial derivatives of $f$. In our formulation, both the  Euler-Maclaurin formula and  the remainder are given as a sum  over  \emph{all} faces of the polytope (not only over vertices) of integrals over cones with apex the affine spaces generated by those faces (see Theorem \ref{thm:4}). 
Moreover, since our formula generalizes to symbols (in the sense of Hormander \cite{H}), we show that, in the case of
polynomial functions, we also obtain an exact Euler-Maclaurin formula for the sum of the values of a polynomial function $p$ over the integer points of the polytope. This relation is a weighted version of the exact Euler Maclaurin formulas obtained in \cite{BrV} and in \cite{KSW}.
\section{Critical set of the function $\vert\vert\mu_{\varepsilon} \vert\vert^2$}
\label{se:cset} Let $(M,\omega)$ be a compact connected symplectic
manifold equipped with a Hamiltonian action of a torus $T$. Denoting
by $\mathfrak{t}$ the Lie algebra of $T$ and by $\mathfrak{t}^*$ its dual
space, we consider the moment map $\mu:M\to \mathfrak{t}^*$ associated
to this action (that is, the $T$-equivariant map  determined, up to
a constant, by the equation $d \langle \mu, X\rangle=\iota(X_M)
\omega$\, for all $X\in \mathfrak{t}$). Since the action of $T$ on
$\mathfrak{t}$ is trivial, the perturbed map $\mu_\varepsilon := \mu
-\varepsilon$ for $\varepsilon \in \mathfrak{t}^*$ is also a moment map
for the action of $T$.

Hereafter we will choose a scalar product on $\mathfrak{t}^*$ (which
induces a linear isomorphism $j: \mathfrak{t} \to \mathfrak{t}^*$) and
consider two kinds of \emph{orthogonality} on $\mathfrak{t}^*$: the
orthogonality resulting from the duality between $\mathfrak{t}$ and
$\mathfrak{t}^*$ and the $j$-orthogonality defined by the scalar
product. We will begin by reviewing the structure of the critical
points of the moment map and then define an index set $\PB$ which
will enable us to define a partition of $\mathrm{Cr}(\vert\vert
\mu_\varepsilon\vert\vert^2)$, the set of critical points  of the
function $\vert\vert \mu_\varepsilon\vert\vert^2$.

Following Paradan in \cite{P}, we will consider a slightly modified
definition of critical point. By the usual definition, a point $x$
is a critical point of the moment map $\mu$ if its stabilizer,
$\mathrm{Stab}(x)$, contains a subtorus of dimension $1$. However,
if the action is not effective, all points of $M$ will be critical
points of $\mu$. To avoid this situation we take the subgroup
$S_M:=\cap_{x\in M} \, \mathrm{Stab}(x)$, called the \emph{generic
stabilizer}, and make the following definition:
\begin{Definition}
The critical points of the moment map $\mu:M\to \mathfrak{t}^*$ are the
points  $x\in M$ for which $\mathrm{Stab}(x)/S_M$ is not finite.
\end{Definition}
Let $T^\prime$ be a subtorus of $T$ containing $S_M$ such that
$T^\prime/S_M$ is not finite. Every connected component $F^\prime$
of $M^{T^\prime}$ is a symplectic submanifold of $M$ and
$P:=\mu(F^\prime)$ is a convex polytope in $\mathfrak{t}^*$ equal to the
convex hull of the image of the fixed points of $T$ contained in
$F^\prime$ (c.f. \cite{At} and \cite{GS}). Moreover, the Lie algebra
of $T^\prime$ is included in $P^\perp$ (the set of vectors of
$\mathfrak{t}$ orthogonal to $P$), and, denoting by $T_P$ the subtorus
of $T$ generated by $\mathrm{exp}(P^\perp)$,  we have that the Lie
algebra $\mathfrak{t}_P$ of $T_P$ is equal to $P^\perp$. The manifold
$F^\prime$ is then a connected component of $M^{T_P}$ where $T/T_P$
acts quasi-freely, that is, the generic stabilizer of $F^\prime$ is
a subgroup of $T$ and its identity component is equal to $T_P$.
Knowing this and denoting by  $\mathrm{Aff}(P)$ the affine subspace of $\mathfrak{t}^*$
generated by $P$, we consider the following sets:
\begin{align}
& \label{set:1} \mathcal{B^\prime}:= \{ \text{convex polytopes} \, P \subset \mathfrak{t}^* \,
\text{for which there exists a connected}\\
& \nonumber \hspace{1.5cm}  \text{component}\,\,  F^\prime \,\,\, \text{of} \,\,\, M^{T_P}\,
\, \text{with} \,\,\, \mu(F^\prime)=P \};
\end{align}
\begin {align}
& \label{set:2}  \PB:=  \{ \mathrm{Aff}(P) \mid P \in \PB^\prime \}; \\
&\label{set:3} \PB^\prime_{\Delta}:=\{ P\in \PB^\prime\mid P\subset
\Delta \, \,\, \text{and} \,\,\, \mathrm{dim}(P) <
\mathrm{dim}(\Delta) \}, \,\, \text{for}\,\,\,
\Delta\in \PB; \\
& \label{set:4} \mathfrak{t}^*_\mathrm{reg}:=  \mathfrak{t}^* \setminus
\bigcup_{P \in \PB^\prime \setminus \mu(M)} \, P,
\,\,\,\,\,\, \,\,\,\,\,\, \,\,\,\,\,\, \Delta_\mathrm{reg}:=  \Delta
\setminus
\bigcup_{P \in \PB^\prime_\Delta}\, P;
\end{align}
\begin{align}
&  \label{set:5} W_\Delta:= \Delta_\text{reg}  + j(\mathfrak{t}_\Delta),\,\,\, \text{for}\,\,
\, \Delta\in \PB;\\
&  \label{set:6} W:= \bigcap_{\Delta\in \PB} \, W_\Delta\,\, .
\end{align}
Note that the set $\PB^\prime$ contains all the
faces of the polytope $\mu(M)$ and that, if $P$ is a polytope in
$\PB^\prime$, all its faces are also in $\PB^\prime$. Again, for
$\Delta \in \PB$, we will denote by $T_\Delta$ the subtorus of $T$
generated by $\mathrm{exp}(\Delta^\perp)$. With this notation we
have the following proposition which characterizes the critical set
of $\vert\vert \mu_\varepsilon\vert\vert^2$ (c.f. \cite{K} and
\cite{P} for details):
\begin{Proposition}\emph{(Kirwan)}
For every $\varepsilon \in \mathfrak{t}^*$,
$$
\mathrm{Cr}(\vert\vert \mu_\varepsilon\vert\vert^2)= \bigcup_{\Delta
\in \PB} M^{T_\Delta} \cap \mu^{-1}(\beta(\varepsilon,\Delta)),
$$
where, for an affine subspace $\Delta$ of $\mathfrak{t}^*$,
$\beta(\varepsilon,\Delta)$ is the orthogonal projection of
$\varepsilon$ on $\Delta$. Moreover, for every $\varepsilon \in
W_\Delta$, the set
\begin{equation}
\label{set:7}
C_\Delta^\varepsilon:= M^{T_\Delta} \cap \mu^{-1}(\beta(\varepsilon,\Delta))
\end{equation}
is a submanifold of $M$ on which $T/T_\Delta$ acts locally freely,
the set $W$ is dense in $\mathfrak{t}^*$ and, for every $\varepsilon \in
W$, the submanifolds $C_\Delta^\varepsilon$, for $\Delta \in \PB$,
form a partition of $\mathrm{Cr}(\vert\vert
\mu_\varepsilon\vert\vert^2)$.
\end{Proposition}
Since the group $T/T_\Delta$ acts locally freely on the manifold
$C_\Delta^\varepsilon$, we can define the quotient
$M^\varepsilon_\Delta:= C_\Delta^\varepsilon/(T/T_\Delta)$  which
will be an orbifold. Moreover, for each connected component $F$ of
$C_\Delta^\varepsilon$, we will consider the subgroup
$S^{\Delta}(F):= \cap_{x\in F}\, \mathrm{Stab}(x)$ and  the map
$F\mapsto \vert S^{\Delta}(F)\vert$ (defining a locally constant
function on $M^\varepsilon_\Delta$) which will be called $\vert
S^\Delta\vert$.
\begin{Remark}
\label{rmk:toric} Note that, if $\mathrm{dim}(T)=1/2\mathrm{dim}(M)$
(i.e. if $M$ is a toric manifold), then
$M^\varepsilon_\Delta=\mu^{-1}(\beta(\varepsilon,\Delta))/T$ is either empty or a
point.
\end{Remark}
From now on we will consider $\varepsilon \in W$.

\section{Localization Formulas}\label{sec:3}
\subsection{Equivariant cohomology}
Let us begin by reviewing the different $T$-equivariant de Rham
complexes on $M$. We have three spaces of equivariant differential
forms on $M$,
$\Omega_T^*(M)\subset \Omega_T^{\infty}(M) \subset \Omega_T^{-\infty}(M)$,
respectively with polynomial, smooth and generalized coefficients.
The model $\Omega_T^*(M)$ is due to Cartan while the other two were
studied by Berline, Duflo, Kumar and Vergne (see \cite{BGV},
\cite{BV1}, \cite{DV} and \cite{KV}). On $\Omega_T^\infty(M)$ we
have the differential $d_T$ defined by
$$
(d_T \alpha) (X) := (d - \iota(X_M))(\alpha(X))
$$
for every $\alpha\in \Omega_T^\infty(M)$ and $X\in \mathfrak{t}$. The
corresponding de Rham $T$-equivariant complexes on $M$ are defined
as $H_T^\infty(M):=(\Omega_T^\infty(M),d_T)$  and $H_T^*(M):=(\Omega_T^*(M),d_T)$,
and called the $T$-equivariant cohomology with $\Cinf$ and
polynomial coefficients. Moreover, the differential $d_T$ defined on
$\Omega_T^\infty(M)$ extends to $\Omega_T^{-\infty}(M)$. Indeed, if
$\{e^1,\ldots, e^r \}$ is a basis for $\mathfrak{t}$ then, for every
$\eta\in  \Omega_T^{\infty}(M)$, we have
$$
\langle d_T(\eta),\phi \rangle := d \langle \eta,\phi \rangle -
\sum_{k=1}^{r}\iota(e_M^k) \langle \eta,\phi \cdot X_k \rangle,
$$
for every compactly supported density $\phi$  on $\mathfrak{t}$, where
$X_1,\ldots, X_r$ are the coordinate functions on $\mathfrak{t}$. The
corresponding de Rham $T$-equivariant complex is denoted
$H_T^{-\infty}(M):=(\Omega_T^{-\infty(M)},d_T)$ and called the
$T$-equivariant cohomolo\-gy on $M$ with generalized
coefficients. We can also consider equivariant classes of compact
support obtaining $T$-equivariant compact support de Rham complexes:
$H_{T,\mathrm{cpt}}^*(M):=(\Omega_{T,\mathrm{cpt}}^*(M),d_T)$,
$H_{T,\mathrm{cpt}}^\infty(M):=(\Omega_{T,\mathrm{cpt}}^\infty(M),d_T)$
and
$H_{T,\mathrm{cpt}}^{-\infty}(M):=(\Omega_{T,\mathrm{cpt}}^{-\infty}(M),
d_T)$.
\subsection{Equivariant Euler classes}
Let $p:E\to M$ be an oriented $T$-bundle and let
$p_*:\Omega_{\mathrm{cpt}}(E) \to \Omega(M)$ be integration along
the fibers. There is a unique equivariant class $\mathit{u}\in
H_{T,\mathrm{cpt}}^\infty(E)$ such that $p_*\mathit{u}=1$ on $M$
called the equivariant Thom class of $E$ and denoted by
$\mathrm{Thom}_T(E)$. Its restriction to $M$ is the equivariant
Euler class of the bundle $E$, i.e.
$
e_T(E):=i^*\mathrm{Thom}_T(E)
$,
where $i:M\to E$ is the inclusion map.

Let us assume now that there is an element $\beta\in \mathfrak{t}$ for
which the zero set of the vector field on $E$ generated by $\beta$
is equal to $M$ and let $T_\beta$ be its stabilizer in $T$.  Then
there is a $T_\beta$-equivariant class with generalized coefficients
$e_\beta^{-1}(E)$ such that $e_\beta^{-1}(E)e_T(E)=1$ (see \cite{P}
for details).
\begin{Example}
Consider the trivial bundle $E:=M\times \C$ equipped with a
$T$-action which is trivial on $M$ and which, on $\C$, is determined by the
weight $\alpha\in \mathfrak{t}^*$ (that is, $\mathrm{exp}(X) \cdot z:= e^{i
\langle \alpha,X \rangle} z$, for $X\in \mathfrak{t}$ and $z\in \C$).
Choosing $\beta\in \mathfrak{t}$ such that $\langle \alpha,\beta \rangle
\neq 0$, we have
$
e_T(E)(X)=-\frac{1}{2\pi} \langle \alpha,X \rangle
$
and, ``polarizing'', that is, taking $\alpha^+:=\epsilon_\beta \,
\alpha$ with $\langle \alpha^+,\beta \rangle > 0$ and
$\epsilon_\beta=\pm 1$, we have
$$
e_\beta^{-1}(E)(X) = 2\pi i \,\epsilon_\beta \int_0^\infty e^{i
\langle \alpha^+,X \rangle\, t} \, dt
$$
as generalized functions. Taking the Fourier transform we obtain the equality of measures on $\mathfrak{t}^*$,
$\mathcal{F}(e_\beta^{-1}(E))=2\pi i\, \epsilon_\beta \, H_{\alpha^+}$,
where $H_{\alpha^+}$ is the Heaviside measure associated to $\alpha^+$ defined by
$$
\langle H_{\alpha^+}, \phi \rangle = \int_0^\infty \phi(u \alpha^+) \, du,
$$
for every $\phi$ in the Schwartz space of rapidly decreasing functions on $M$.
\end{Example}
\begin{Example}
If $M=\{F\}$ is a single point, fixed by the action of $T$, the
bundle $E$ decomposes as a sum of non-trivial $2$-dimensional real
representations of $T$,
$E:=L_1 \oplus \cdots \oplus L_k \to F$,
with the action of $T$  on each $L_j$ determined by a weight
$\alpha_j\in \mathfrak{t}^*$. Following Paradan (cf. Proposition 4.8 in
\cite{P}) we obtain the expression for the Fourier transform of
$e_\beta^{-1}(E)$:
$$
\mathcal{F}(e_\beta^{-1}(E))=(2\pi i )^k \, \epsilon_\beta \,
H_{\alpha_1^+} * \cdots * H_{\alpha_k^+},
$$
where we polarize each $\alpha_j$ according to some $\beta\in \mathfrak{t}$
(such that $\langle \alpha_j, \beta \rangle \neq 0$ for
$j=1,\ldots,k$), obtaining $\alpha_j^+:=\epsilon_\beta^j \,
\alpha_j$ with $\epsilon_\beta^j=\pm 1$, and we take
$\epsilon_\beta:= \prod_{j=1}^k \epsilon_\beta^j$. Note that $*$
denotes the convolution product. This measure, supported on the cone
$\R^+ \alpha_1^+ + \cdots + \R^+ \alpha_k^+$, is defined by
$$
\langle H_{\alpha_1^+} * \cdots * H_{\alpha_k^+}, \phi \rangle =
\int_0^\infty \cdots   \int_0^\infty \phi(\sum_{i=1}^k u_i
\alpha_i^+) \, du_1 \ldots du_k,
$$
for every rapidly decreasing function on $M$.
\end{Example}
\subsection{Localization}
Using the sets $\PB$, $W$ and $C_\Delta^\varepsilon$ defined in
(\ref{set:2}), (\ref{set:6}) and  (\ref{set:7}) of
Section~\ref{se:cset}, and the orbifold $M_\Delta^\varepsilon =
C_\Delta^\varepsilon/(T/T_\Delta)$,  Paradan proves the following
localization theorem:

\begin{Theorem}\emph{(Paradan)}
Let $\varepsilon \in W$ and let $\eta \in \Omega_T^\infty(M)$ be a
closed form. Then, on $C^{-\infty}(\mathfrak{t})$ we have
\begin{equation*}
\int_M \eta = \sum_{\Delta \in \PB}\, I_\Delta^\varepsilon(\eta),
\end{equation*}
where $I_\Delta^\varepsilon(\eta)$ is the generalized function supported on
$\mathfrak{t}_\Delta$ defined by
\begin{align}
\label{eq:paradanloc}
\begin{split}
 I_\Delta^\varepsilon(\eta)(X_1&+X_2) = \\
 &=(2\pi i)^{\mathrm{dim} \Delta}
 \int_{M_\Delta^\varepsilon} \frac{1}{\vert S^\Delta \vert} k_\Delta(\eta)(X_1)
 e_{\beta_\Delta}^{-1}(E_\Delta)(X_1) \diamond \delta(X_2-w_\Delta),
\end{split}
\end{align}
where
\begin{enumerate}
\item[(i)] the variables $X_1, X_2$ are respectively in $\mathfrak{t}_\Delta$ and
$\mathfrak{t}/\mathfrak{t}_\Delta$ (note that, for each $\Delta \in \PB$,
$\mathfrak{t}$ decomposes as a sum of vector spaces $t_\Delta$ and
$\mathfrak{t}/\mathfrak{t}_\Delta$, where $\mathfrak{t}_\Delta$ is the Lie
algebra of the subtorus $T_\Delta$ generated by
$\mathrm{exp}(\Delta^\perp)$);
\item[(ii)] $k_\Delta: H_T^\infty(M) \to H_{T_\Delta}^\infty(M_\Delta^\varepsilon)$
is the Kirwan map (see \cite{K});
\item[(iii)]$\beta_\Delta:=j^{-1}(\beta(\varepsilon,\Delta) - \varepsilon)$,
where $\beta(\varepsilon,\Delta)$ is the orthogonal projection of
$\varepsilon$ on $\Delta$;
\item[(iv)] $E_\Delta:= N_\Delta/(T/T_\Delta)$, where  $N_\Delta$ is the normal
bundle of $M^{T_\Delta}$ inside $M$, restricted to $C_\Delta^\varepsilon$;
\item[(v)] the operator $\diamond$:
\begin{eqnarray*}
\Omega_{T_\Delta}^{-\infty}(C_\Delta^\varepsilon) \times   \Omega_{T/T_\Delta}^{-\infty}
(C_\Delta^\varepsilon) &  \to & \Omega_T^{-\infty}(C_\Delta^\varepsilon) \\
     (\eta,\nu)                                        &  \mapsto & \eta \diamond \nu
\end{eqnarray*}
is defined by
$$
\langle \eta \diamond \nu, \phi(X)dX \rangle := \langle \eta ,
\langle \nu, \phi(X_1 + X_2) dX_1 \rangle dX_2 \rangle
$$
for every density $\phi(X)$ of compact support on $\mathfrak{t}$;
\item[(vi)] $w_\Delta$ is the equivariant curvature of the principal orbibundle
$C_\Delta^\varepsilon \to M_\Delta^\varepsilon$;
\item[(vii)] the equivariant form $\delta(X_2-w_\Delta) \in
\Omega_{T/T_\Delta}^{-\infty}(M_\Delta^\varepsilon)$ is defined by
$$
\langle \delta(X_2-w_\Delta), \phi(X_2) d X_2 \rangle = \phi(w_\Delta)
\mathrm{vol}(T/T_\Delta, dX_2),
$$
for every function $\phi \in \Cinf(\mathfrak{t}/\mathfrak{t}_\Delta)$, where
$\mathrm{vol}(T/T_\Delta, dX_2)$ is the volume of the group
$T/T_\Delta$ with respect to the Haar measure compatible with
$dX_2$.
\end{enumerate}
\end{Theorem}
\begin{Example}
\label{ex:1}
If $\Delta=\{p \}$ is a vertex of the polytope $\mu(M)$,
then $C_\Delta^\varepsilon=\mu^{-1}(p)$ is a connected component $F$ of
$M^T$ and
$$
I_{\{p\}}^\varepsilon(\eta)(X)= \int_F i^*_F(\eta)(X) e_{\beta_p}^{-1}(N_F) (X),
$$
where $N_F$ is the normal bundle of $F$ inside $M$ and
$\beta_p=j^{-1}(p - \varepsilon)$. If, in addition, the action of
$T$ is toric, then $F$ is an isolated point. Moreover, taking
$\eta=e^{i\omega^\sharp}$, where $\omega^\sharp$ is the equivariant
symplectic form on $M$ defined by $\omega^\sharp (X):=\omega(X) +
\langle \mu, X \rangle$, we obtain
\begin{equation}
 I_{ \{ p\} }^\varepsilon (e^{i\omega^\sharp})(X)   =   \epsilon_p \,
 e^{i\langle p,X \rangle}\, \prod_{j=1}^{\mathrm{dim}(M)/2} \int_0^\infty
 e^{i\langle \alpha_j^+,X \rangle t} \, dt,
\end{equation}
where the $\alpha_j^+$'s are the polarized weights of the action of
$T$ on the normal bundle of $F$ (i.e. the polarized edge vectors at
$p$) and $\epsilon_p$ is the sign obtained by polarization. Taking its
Fourier transform we get
\begin{equation}
\label{eq:2} \mathcal{F} (I_{ \{ p \} }^\varepsilon
(e^{i\omega^\sharp})) =  \epsilon_p \, \delta_p * H_{\alpha_1}^+ *
\cdots * H_{\alpha_n}^+,
\end{equation}
where $\delta_p$ is the Dirac measure on $p \in \mathfrak{t}^*$.
Moreover, the measure  (\ref{eq:2}) is supported on the polarized
cone
$
\bfC_{\{p\}}^\sharp:=p + \R^+ \alpha_1^+ + \cdots + \R^+ \alpha_n^+
$.
\end{Example}
\begin{Example}\label{ex:dec}
For a toric manifold $M^{2n}$ with moment map $\mu$ and for  $\eta=
e^{i\omega^\sharp}$, where again $\omega^\sharp$ is the equivariant
symplectic form on $M$ (cf. Example~\ref{ex:1}), the reduced space
$M_\Delta^\varepsilon$ for $\varepsilon
\in \mathfrak{t}^*$ is either empty or a single point. Moreover, denoting by $\beta_1(\varepsilon,
\Delta)$ the orthogonal projection of  $\beta(\varepsilon, \Delta)$
onto  $\mathfrak{t}_\Delta^*$, formula \eqref{eq:paradanloc} becomes
\begin{align}
\label{eq:3}
\begin{split}
I_\Delta^\varepsilon (X_1 &+ X_2) \\
&= (2\pi i)^{\dim \Delta}
\epsilon_\Delta \left( \frac{e^{i \langle \beta_1(\varepsilon,
\Delta), X_1 \rangle}}{\vert S^\Delta\vert} \prod_{j=1}^{r_\Delta}
 \int_0^\infty e^{i\langle \alpha_{\Delta,j}^+,X_1 \rangle t}\,
dt \right) \diamond \, \delta_0(X_2),
\end{split}
\end{align}
whenever $\beta(\varepsilon,\Delta)\cap \mu(M)$ is nonempty, where
$\epsilon_\Delta =\prod_{j=1}^{r_\Delta} \epsilon_\Delta^j$ with
$\alpha_{\Delta,j}^+=\epsilon_{\Delta}^j \alpha_{\Delta,j}$ is the
sign obtained by polarization,  where $r_\Delta$ is the codimension of
$\Delta$, where $\vert S^\Delta\vert$ is the order of the orbifold
structure group of the point in $M_\Delta^\varepsilon$ inside the
toric orbifold $E_\Delta:=N_\Delta/T/T_\Delta$, and where the
$\alpha_{\Delta,j}^+$'s are the polarized weights of the action of
$T_\Delta$ on $E_\Delta$ restricted to the normal orbibundle of the
fixed point. Equivalently, $\vert S^\Delta\vert$ is the order of the
orbifold structure group of the point in $M_\Delta^\varepsilon$
inside the reduced space $\mu_{T/T_\Delta}^{-1}
(\beta_2(\varepsilon, \Delta))/T/T_\Delta$, where $\mu_{T/T_\Delta}$
is the moment map for the $T/T_\Delta$-action on $M$, and
$\beta_2(\varepsilon, \Delta)$ is the orthogonal projection of
$\beta(\varepsilon, \Delta)$ onto $(\mathfrak{t}/\mathfrak{t}_\Delta)^*$.

Taking the  Fourier transform of (\ref{eq:3}) we obtain
\begin{equation}
\label{eq:4}
(2\pi i)^n \frac{\epsilon_\Delta}{\vert S^\Delta\vert}
\left( (\delta_{\beta_1(\varepsilon, \Delta)} *
H_{\alpha_{\Delta,1}^+} * \cdots *  H_{\alpha_{\Delta,r}^+})
\diamond \, \mathbf{1}_{(\mathfrak{t}/\mathfrak{t}_\Delta)^*} \right),
\end{equation}
which is supported on the polarized cone
$
\bfC_{\beta_1(\varepsilon, \Delta)}^\sharp:= \beta_1(\varepsilon,
\Delta) + \R^+ \alpha_{\Delta,1}^+ +\cdots +
\R^+\alpha_{\Delta,r_\Delta}^+
$.
Moreover, changing variables, we obtain
\begin{align*}
\langle  H_{\alpha_{\Delta , 1}^+} * \cdots * H_{\alpha_{\Delta ,
r_\Delta}^+} , \phi \rangle & = \int_0^\infty \cdots   \int_0^\infty
\phi(\sum_{i=1}^{r_\Delta} u_i \, \alpha_{\Delta,i}^+) \, du_1
\ldots du_{r_\Delta}= \\ & = \frac{1}{\vert
\det{(\alpha_{\Delta,i}^+)_i}\vert } \int_{\bfC_0^\sharp} \phi =
\frac{1}{\vert \det{(\alpha_{\Delta,i}^+)_i}\vert } \langle
\bfone_{\bfC_0^\sharp} , \phi \rangle,
\end{align*}
for any rapidly decreasing function $\phi$, where $\bfone_{\bfC_0^\sharp}$
is the characteristic function of the cone $\bfC_0^\sharp := \R^+
\alpha_{\Delta,1}^+ + \cdots + \R^+\alpha_{\Delta,r_\Delta}^+$.
However, since the $\alpha_{\Delta , i}$'s are the weights of the
action of $T_\Delta$ on the toric orbifold $E_\Delta$ at the fixed
point, we have
$\vert \det{(\alpha_{\Delta,i}^+)_i}\vert =\frac{1}{\vert S^\Delta \vert}$
and so \eqref{eq:4} becomes $(2\pi i)^n  \epsilon_\Delta \,
\bfone_{\bfC_0^\sharp}  \diamond \,
\mathbf{1}_{(\mathfrak{t}/\mathfrak{t}_\Delta)^*}.$ Indeed, denoting by
$\hat{\ell}^*$ the lattice dual to the weight lattice $\hat{\ell}$
of $T_\Delta$, the orbifold structure group $\Gamma$ of the fixed
point is isomorphic to $\ell/\hat{\ell}$, where $\ell$ is the
lattice of circle subgroups of $T_\Delta$ (cf. \cite{LT} for
details), implying that 
$$
\vert S^\Delta \vert = \vert \Gamma \vert = \vert
\det{\hat{\ell}} \vert = \frac{1}{\vert \det{\ell}
\vert} = \frac{1}{\vert \det{(\alpha_{\Delta,i}^+)_i \vert}}
$$
(see for example \cite{C}). On the other hand, since the Fourier transform of
$\int_M e^{i \omega^\sharp t}$ is the direct image $\mu_*(dm_L)$ of
the Liouville measure $dm_L:= \omega^n/n!$ on $M$ (which is
supported on $\mu(M)$), we obtain (modulo $(2\pi i)^n$)
\begin{equation}
\label{eq:paradandecomp} \bfone_{\mu(M)}= \sum_{\Delta\in \PB}
\varphi(\varepsilon,\Delta) \, \epsilon_\Delta \,
\bfone_{\bfC_\Delta^\sharp},
\end{equation}
where $\varphi(\varepsilon,\Delta)$ is equal to $1$ when
$\beta(\varepsilon,\Delta)\cap \mu(M)$ is nonempty, and zero
otherwise.
Note that Formula \eqref{eq:paradandecomp} is valid up to boundary
effects on the polytope $\mu(M)$.
\end{Example}
\section{Polytope decompositions}
\label{se:decompositions} In this section we will show that the
polytope decomposition for
Delzant polytopes that was obtained in \eqref{eq:paradandecomp}, remains valid for any compact convex
simple  polytope. Moreover, we will give a
weighted version of this decomposition that also holds on the boundary of the
polytope.

\begin{figure}[h]
  \centering
 \includegraphics[scale=.45]{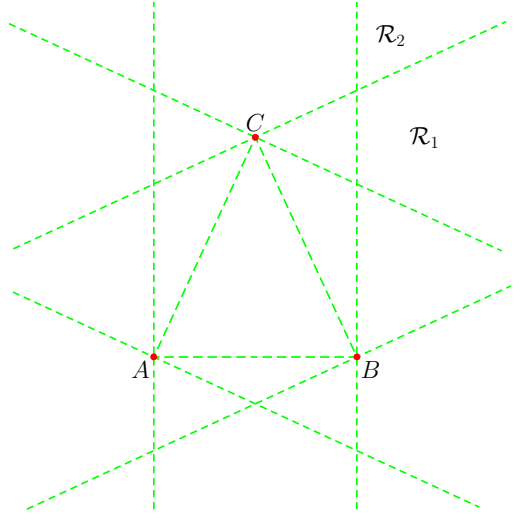}\\
  \caption{Paradan regions for a triangle.}\label{fi:acosta0}
\end{figure}

Hereafter, we will consider the usual Euclidean inner product
$\langle\,,\rangle$ of $\R^d$. Let $P$ be a compact convex simple
polytope in $\R^d$ and let $\PB^\prime$ be the set of faces of $P$.
For each $F\in \PB^\prime$  we write $\Delta_F$ for the affine
subspace of $\R^d$ generated by $F$. Then, just as in
Section~\ref{se:cset}, we have the following sets:
\begin{align}
& \label{set:3.2}  \PB:=  \{ \Delta_F \mid F \in \PB^\prime \}; \\
& \label{set:3.3}  \PB^\prime_{\Delta}:=\{ F\in \PB^\prime\mid F\subset
\Delta \,\,\, \text{and} \,\,\, \mathrm{dim}(F) < \mathrm{dim}(\Delta) \},
\,\, \text{for}\,\,\,\Delta\in \PB; \\
& \label{set:3.4}   \Delta_\mathrm{reg}:=  \Delta
\setminus \bigcup_{F \in \PB^\prime_\Delta}\, F;\\
&  \label{set:3.5} W_\Delta:= \Delta_\mathrm{reg}  + \Delta^\perp\,\,\,
\text{for}\,\,\, \Delta\in \PB;
\end{align}
\begin{align}
&  \label{set:3.6} W:= \bigcap_{\Delta\in \PB} \, W_\Delta.
\end{align}
The set $W$ is a disjoint union of open sets which  we will call
\emph{Paradan regions} (see Figure~\ref{fi:acosta0} for an
illustration). In fact, $W$ is the complement in $\R^d$ of a finite
set of walls of codimension $1$, $W^c = E_1\cup
\cdots \cup E_K \cup \{\text{facets of P}\}$, where each wall $E_i$ is contained in
a hyperplane perpendicular to a family of elements of $\PB$. Note that, in the
case of a moment polytope $\mu(M)$, the set $W$ is the same as in
(\ref{set:7}).

\subsection{Tangent cones}
\begin{figure}[h]
  \centering
  \includegraphics[scale=.45]{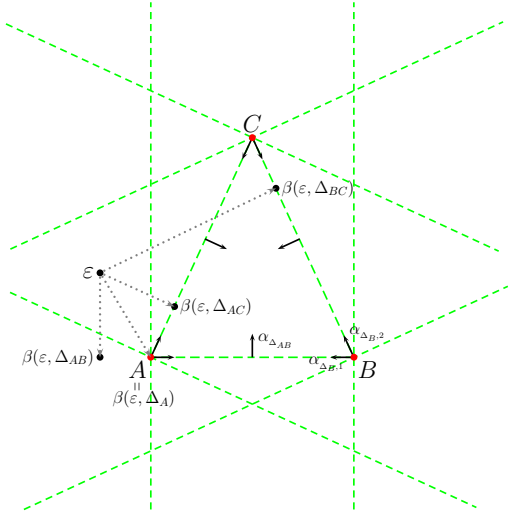}\\
  \caption{Projections and generating vectors for some faces of a triangle.}
  \label{fi:acosta1}
\end{figure}

For each $\Delta \in \PB$, we define the tangent cone of $P$ at $\Delta$ by
\begin{equation*}
\bfC_{\Delta}:=\{y+r(x-y)\,|\, r\ge0,\, y\in \Delta,\, x\in P\}.
\end{equation*}
It is a full dimensional cone with apex $\Delta$ (i.e. $\Delta$ is
the maximal affine space contained in $\bfC_\Delta$). Taking
$\varepsilon$ in $\R^d$, we denote by $\beta(\varepsilon,\Delta)$
its orthogonal projection onto the affine subspace $\Delta$,  and we
take the intersection of $\beta(\varepsilon,\Delta) + \Delta^\perp$
(the orthogonal space of $\Delta$ at  $\beta(\varepsilon,\Delta)$)
with the tangent cone of $\Delta$,
\begin{equation*}
\bfC_{\Delta^\perp,\varepsilon}:=(\beta(\varepsilon,\Delta) +
\Delta^\perp) \cap \bfC_\Delta,
\end{equation*}
which is now a pointed cone\footnote{A cone with a single point as apex (also called vertex).} with
vertex $\beta(\varepsilon,\Delta)$. The cone $\bfC_\Delta$ is the
direct product of the affine space $\Delta$ and the pointed cone
$\bfC_{\Delta^\perp,\varepsilon}$. Then, considering vectors
$\alpha_{\Delta,j}\in \R^d$ along the edges of $\bfC_{\Delta^\perp,\varepsilon}$, pointing
away from the vertex ($j=1,\ldots,r_\Delta$, where $r_\Delta =
\dim{\bfC_{\Delta^\perp,\varepsilon}}= \text{codim}\, \Delta$), the
tangent cone $\bfC_\Delta$ can be written as
\begin{equation}\label{eq:conealongDelta}
\bfC_{\Delta}= \Delta + \bfC_{\Delta^\perp,\varepsilon} =\Delta +
\sum_{j=1}^{r_\Delta}\R^+\alpha_{\Delta,j},
\end{equation}
that is, $\bfC_\Delta$ is the cone along $\Delta$ which contains $P$
and is  bounded by the affine spaces in $\PB$ 
which contain  $\Delta$. The vectors $\alpha_{\Delta,j}$ (which are
only determined up to a positive scalar) will be called the
\emph{generators} of $\bfC_\Delta$ (see Figure~\ref{fi:acosta1} for
an illustration).
\bigskip

\subsection{Polarization} We will ``polarize'' the tangent cones
in the following way:  let us consider a point  $\varepsilon \in W$; as
a direct consequence of the definition of $W$, for each $\Delta \in
\PB$ we have $\langle \beta_\Delta,  \alpha_{\Delta,j}\rangle\neq 0$ for all
$j=1,\ldots,r_\Delta$, where
$\beta_\Delta:=\beta(\varepsilon,\Delta) - \varepsilon$; polarizing
the vectors $\alpha_{\Delta,j}$ according to $\beta_\Delta$, that
is, taking the vectors $\alpha^+_{\Delta,j}
=\epsilon_{\beta_\Delta}^j \alpha_{\Delta,j}$ with $\langle
\alpha_{\Delta,j}^+ , \beta_{\Delta}\rangle > 0$  and
$\epsilon_{\beta_\Delta}^j =\pm 1$, we define
$\bfC^\sharp_{\Delta}$, the \emph{polarized tangent cone} of $P$ at
$\Delta$ by
\begin{equation}\label{eq:polarizedcone}
\bfC^\sharp_{\Delta}:=\Delta + \sum_{j=1}^{r_\Delta} \R^+ \alpha^+_{\Delta,j}.
\end{equation}

\subsection{Weighted characteristic functions} Let us now see how to assign weights to each affine space in $\PB$
in order to obtain a weighted version of \eqref{eq:paradandecomp}
that also holds for points in $\partial P$. Take $\Delta_1,\ldots,
\Delta_N$, the codimension-$1$ elements of $\PB$, that is, the
affine subspaces generated by the facets of $P$. Each $\Delta\in
\PB$ that is generated by a non-trivial face of $P$ (i.e. such that
$\Delta \neq \varnothing, \R^d$) can be described as an intersection
$$
\bigcap_{i\in J_\Delta} \,\Delta_i,
$$
where $J_\Delta$ denotes the index set of the hyperplanes $\Delta_i$
that contain $\Delta$. Note that, since $P$ is simple, the number of
elements of $J_\Delta$ is equal to $r_\Delta$, the codimension of
$\Delta$. To each $\Delta_i$ we assign an arbitrary complex number
$q_i$ and to each affine space $\Delta\in \PB\setminus \R^d$ we
assign the value $\prod_{i\in J_\Delta} q_i$. Moreover, to
$\Delta=\R^d$ we assign the value $1$.  This amounts to defining a
weighted function $w:\R^d \to \C$, given by $w(x)=\prod_{i\in
J_{\Delta_x}}q_i$, where $\Delta_x$ is the smallest-dimension element of
$\PB$ that contains $x$. With this function, we define the
\emph{weighted characteristic function}
\begin{equation}\label{eq:weightedpolytope}
\bfone^w_{P}(x)=\left\{\begin{array}{cc}
                        w(x), & \mbox{if }x\in P \\
                        0, & \mbox{otherwise} \\
                      \end{array}\right. .
\end{equation}

\begin{Remark}\label{le:wcharcone} Similarly,
for each $\Delta\in\PB$, we define  weighted characteristic
functions for the tangent cone $\bfC_\Delta$ and for the polarized
tangent cone $\bfC^\sharp_\Delta$, but now assigning the weight $1-q_i$ to each facet of $\bfC_\Delta$ that is flipped in
the polarization process.
\end{Remark}

The fact that this assignment of weights fits together and makes our
polytope decomposition hold even for points in $\partial P$, relies
heavily on the following Lemma as we will see later in the proof of
Theorem~\ref{thm:1}.
\begin{Lemma}\label{le:qi's} Let $I=\{1,\ldots,n\}$. For any $q_i\in\C$ with $i\in I$, we have
\begin{equation}\label{eq:qi's}
\prod_{i\in I} q_i = 1 + \dsum_{\varnothing\neq J\subset
I}\,\,\prod_{j\in J} (q_j-1).
\end{equation}
\end{Lemma}

\begin{proof} We will use induction on the cardinality $n:=\#(I)$ of $I$. Clearly,
formula \eqref{eq:qi's} is trivial for $n=1$. Let us then assume
that \eqref{eq:qi's} holds for the set $I^\prime=\{1,\ldots,n-1\}$ and show that it holds for 
$I=\{1,\ldots,n\}$. Indeed we have
\begin{align*}
 \dprod_{i\in I} q_i & =  \left( \dprod_{i=1}^{n-1} q_i\right) q_n =
 \left(1 + \dsum_{\varnothing \neq J^\prime\subset I^\prime}
 \,\, \prod_{j\in J^\prime}(q_j-1)\right) (q_n - 1 +1) \\
& =  1 + (q_n - 1)  + \dsum_{\varnothing \neq J^\prime \subset
I'}\,\,\left(\dprod_{j\in J^\prime}(q_j - 1) (q_n - 1) +
\dprod_{i\in J^\prime}(q_i - 1)\right).
\end{align*}
However, on the other hand,
\begin{equation*}
 \dsum_{\varnothing \neq J\subset I}\,\, \dprod_{j\in J} (q_j-1)
=  \dsum_{\substack{J\subset I \\[.5mm] \text{with}\, n\in J}}
  (q_n-1) \dprod_{j\in J\setminus\{n\}}(q_j-1)
 + \dsum_{\varnothing \neq J\subset I'}\,\, \prod_{i\in J}(q_i-1),
\end{equation*}
and the result follows.
\end{proof}

\subsection{Decomposition formulas}
\begin{figure}[h]
  \centering
  \includegraphics[scale=.6]{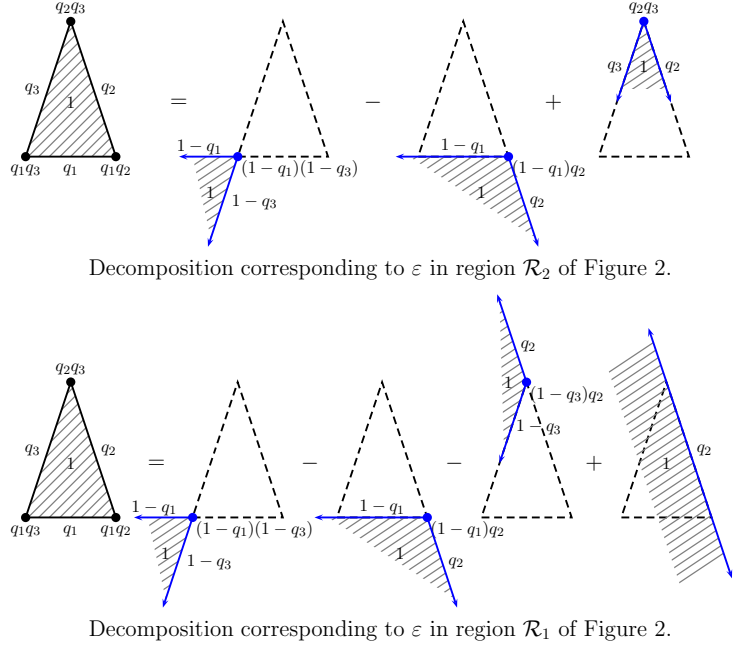}\\
  \caption{Polytope decomposition for a triangle}\label{fi:Paradantriangle}
\end{figure}

Finally, defining $\varphi(\varepsilon,\Delta)$ as
\begin{equation}\label{eq:indicator}
\varphi(\varepsilon,\Delta)=\left\{\begin{array}{ll}
                        1 & \mbox{if }\beta(\varepsilon,\Delta)\cap P
                        \neq \varnothing \\
                        0 & \mbox{otherwise} \\
                      \end{array}\right. ,
\end{equation}
we obtain, for each $\varepsilon\in W$, the following polytope
decomposition formula (See
Figure ~\ref{fi:Paradantriangle} for an illustration). Note that, by
the definition of $\varphi$, this formula only takes into account the
polarized cones $\bfC_\Delta^\sharp$ for which $\beta(\varepsilon,
\Delta) \cap P \neq \varnothing$, that is, those for which the
orthogonal projection of $\varepsilon$ onto $\Delta$ is in $P$.

\begin{Theorem}
\label{thm:1}
For any compact convex simple polytope $P$ of dimension $d$ in $\R^d$ and for any
$\varepsilon$ in $W$, we have
\begin{equation} \label{Paradan}
\bfone^w_{P} = \sum_{\Delta \in \PB} (-1)^{m_\Delta}
\varphi(\varepsilon,\Delta)\bfone^w_{\bfC^\sharp_{\Delta}},
\end{equation}
where the sum is taken over the set $\PB$ of affine spaces generated
by the  faces of $P$, where $\bfC^\sharp_{\Delta}$ is the tangent
cone of $P$ at  $\Delta\in \PB$ polarized with respect to the vector
$\beta(\varepsilon,\Delta)-\varepsilon$ (where
$\beta(\varepsilon,\Delta)$ is the orthogonal projection of
$\varepsilon$ onto $\Delta$), where $m_\Delta$ is the number of
generators of the tangent cone $\bfC_{\Delta}$ whose signs change by
polarization, and where $\bfone^w_{P}$ and
$\bfone^w_{\bfC^\sharp_{\Delta}}$ are the weighted characteristic
functions of the polytope and of the polarized cones respectively.
\end{Theorem}
\begin{proof}
We will prove this formula in two steps. First, for
each $x$ we will find an $\varepsilon$ in a Paradan region for which
formula \eqref{Paradan} holds. Then, we will show that the right
hand side is independent of the choice of $\varepsilon$.

{\bf Step 1:} Suppose $x\notin P$. If $x$ is in $W$ then we choose
$\varepsilon$ on the same Paradan region as $x$. If, on the other
hand,  $x$ is in the complement of $W$, we choose $\varepsilon$ on
any Paradan region contiguous to $x$ inside the complement of $P$.
In both cases, none of the cones $\bfC^\sharp_\Delta$ on the right
hand side of (\ref{Paradan}) contains $x$ and we obtain
$$ 0 =\bfone^w_{P}(x) = \sum_{\Delta \in \PB} (-1)^{m_\Delta}
\varphi(\varepsilon,\Delta)\bfone^w_{\bfC^\sharp_{\Delta}}(x) = 0.$$

If $x$ is on the boundary of $P$, we choose $\varepsilon$ in the
interior of $P$ close enough to $x$ so that
$\varphi(\varepsilon,\Delta)=1$ for every $\Delta \in \PB$ that
contains $x$. In this case, all polarized cones $\bfC^\sharp_\Delta$
that contain $x$ are pointing away from $\varepsilon$. Moreover, all
generators of the corresponding tangent cones $\bfC_{\Delta}$ are
flipped and so $m_{\Delta}=\codim(\Delta)$. Denoting by $\Delta_x$
the smallest dimensional affine subspace in $\PB$ that contains $x$,
we have $\bfone^w_{P}(x)=\prod_{i\in J_{\Delta_x}} q_i$, while the
right hand side of \eqref{Paradan} is
\begin{align*}
& \text{RHS}  = \dsum_{\Delta \in \,\PB_x } (-1)^{m_\Delta}
\varphi(\varepsilon,\Delta) \bfone^w_{\bfC^\sharp_{\Delta}}(x)
= \dsum_{\Delta \in\,\PB_x} (-1)^{m_\Delta} w_\Delta(x) \\ 
& =   1  + \!\!\!\!\!\!\dsum_{\Delta \in\,\PB_x \setminus \R^d}  (-1)^{m_\Delta} \dprod_{i\in J_{\Delta}} (1-q_i)
=  1  + \!\!\!\!\!\!\dsum_{\Delta \in\,\PB_x \setminus \R^d}  \dprod_{i\in J_{\Delta}} (q_i - 1) 
=  1  + \!\!\!\!\! \dsum_{\substack{J \subset J_{\Delta_x} \\ J \neq \varnothing}} \dprod_{j\in \,J} (q_j-1)\,,
\end{align*}

\noindent where $\PB_x\subset \PB$ is the subset of elements of
$\PB$ that contain $x$, where $w_\Delta(x):= \prod_{j\in J_{\Delta}}
(1-q_j)$, and where we used the fact that $m_\Delta = \codim{\Delta}
=\#J_{\Delta}$ and that $\Delta \in \PB_x$ iff $J_\Delta\subset
J_{\Delta_x}$. Consequently, from Lemma~\ref{le:qi's} we conclude
that
\begin{equation*}
\bfone^w_{P}(x) = \sum_{\Delta \in \PB} (-1)^{m_\Delta}
\varphi(\varepsilon,\Delta)\bfone^w_{\bfC^\sharp_{\Delta}}(x).
\end{equation*}

If $x$ is in the interior of $P$ then, taking $\varepsilon=x$, none
of the polarized $\bfC_\Delta^\sharp$ in (\ref{Paradan})  contains $x$
except the one corresponding to the top dimensional face $P$ of the
polytope. In this case, $\varphi(\varepsilon,\R^d)=1$, and, formula
(\ref{Paradan}) evaluated at $x$ gives
$$1= \bfone^w_{P}(x) =\varphi(\varepsilon,\R^d)\bfone^w_{\bfC^\sharp_{\R^d}}(x)
=  1.$$
\begin{figure}[h]
  \centering
 \includegraphics[scale=.60]{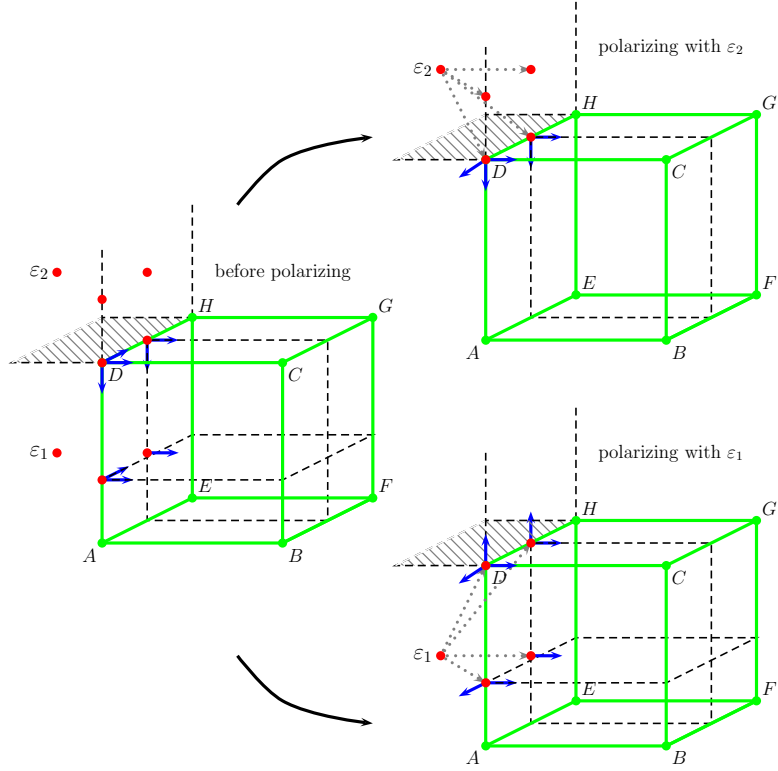}\\
  \caption{Generators of some polarized tangent cones $\bfC_\Delta^\sharp$ for a cube,
  with $\varepsilon_1$ and $\varepsilon_2$ in two contiguous Paradan regions.}
  \label{fi:cubo}
\end{figure}
{\bf Step 2:} Recall that the complement $W^c$ of $W$ is a finite
family of walls of codimension $1$, $W^c= E_1 \cup \cdots \cup E_K
\cup \{\text{facets of P}\}$, where each wall  $E_i$ is perpendicular to a
family $\PB_i$ of affine spaces contained in $\PB$. Let
$\varepsilon_1$ and $\varepsilon_2$ be in two contiguous Paradan
regions $\PR_1$ and $\PR_2$ respectively, and let $E$ be its common
``wall'' (either one of the $E_i$'s or a facet of $P$). Let
$\varepsilon_t$ be any path in $\R^d$ from $\varepsilon_1$ to
$\varepsilon_2$ that crosses a single wall (i.e. $E$) once. When
$\varepsilon_t$ crosses $E$, the sign of
\begin{equation}\label{eq:beta}
\langle\beta(\varepsilon_t,\Delta)-\varepsilon_t,\alpha_{\Delta,k}\rangle
\end{equation}
(for $\Delta\in \PB \setminus \R^d$ and for a generator $\alpha_{\Delta,k}$
of the tangent cone at $\Delta$) flips exactly when $\Delta\cap \partial
P$ is contained in $E$ and $\alpha_{\Delta,k}$ is perpendicular to
$E$. Hence, if $\dim{\Delta}\neq d -1$, this sign flips iff  $\Delta$
is contained in (exactly) one affine space $\widetilde{\Delta}$ of
$\PB_i$ with $\dim \widetilde{\Delta} = \dim \Delta +1$ (see Figure
~\ref{fi:cubo} for an illustration). Indeed, we can take
$$
\widetilde{\Delta}:=\Delta + \text{span}\,{\alpha_{\Delta,k}} \in \PB_i
$$
and unicity follows from dimensional reasons. On the other hand, if
$\dim{\Delta} =d - 1$, the sign of \eqref{eq:beta} flips iff $\Delta
\cap P\subseteq E$. In this case, we define  $\widetilde{\Delta}$ to be the entire
space $\R^d$.
\begin{figure}[h]
   \centering
   \includegraphics[scale=.70]{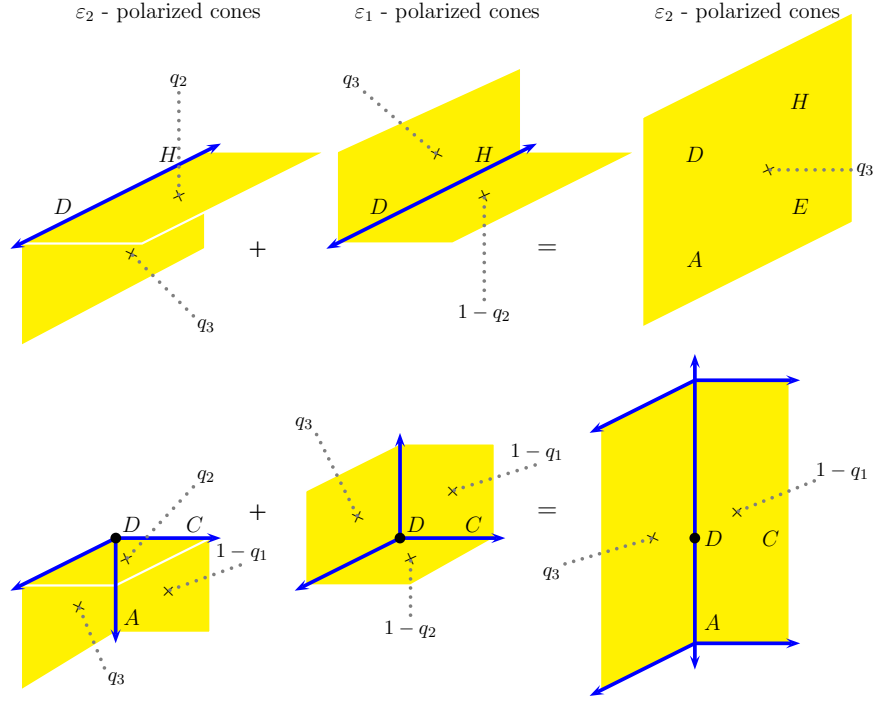}\\
   \caption{Polarized tangent cones $\bfC^\sharp_\Delta$ at the affine subspaces
   $\Delta_D$, $\Delta_{D\!H}$,
   $\Delta_{D\!A}$ and $\Delta_{A\!E\!H\!D}$ for the cube in Figure~\ref{fi:cubo}.}
   \label{fi:cones3d}
 \end{figure}

Let us assume without loss of generality that the sign of \eqref{eq:beta}
flips from negative to positive as $\varepsilon_t$ crosses $E$. In
this case, the polarized tangent cones at $\Delta$ before and after
$\varepsilon_t$ crosses the wall are
$$(\bfC^\sharp_{\Delta})^1 = \Delta + \sum_{j\neq k} \R^+ \alpha^\sharp_{\Delta,j}
- \R^+ \alpha_{\Delta,k} \quad \text{and} \quad
(\bfC^\sharp_\Delta)^2 = \Delta + \sum_{j\neq k} \R^+ \alpha^\sharp_{\Delta,j} +
\R^+ \alpha_{\Delta,k}. $$
Hence, the corresponding contributions
of $\Delta$ to the right hand side of \eqref{Paradan} are
\begin{align*}
\pm \bfone^w_{(\bfC^\sharp_\Delta)^1} &  & \text{and} & & \mp
\bfone^w_{(\bfC^\sharp_\Delta)^2}.
\end{align*}
Note that the union of the two cones $(\bfC^\sharp_\Delta)^1$ and
$(\bfC^\sharp_\Delta)^2$ is the polarized tangent cone at
$\widetilde{\Delta}$, $\bfC^\sharp_{\widetilde{\Delta}}$, for both
$\varepsilon_1$ and $\varepsilon_2$ (cf. Figure~\ref{fi:cones3d}),
and so
$$
\bfone^w_{(\bfC^\sharp_\Delta)^1} +
\bfone^w_{(\bfC^\sharp_\Delta)^2} =
\bfone^w_{\bfC^\sharp_{\widetilde{\Delta}}}.
$$

On the other hand, we have $\beta(\varepsilon_1,\widetilde{\Delta})\cap P \neq
\varnothing$, while $\beta(\varepsilon_2,\widetilde{\Delta})\cap P =
\varnothing$. Hence, the corresponding contributions of
$\widetilde{\Delta}$ to the right hand side of \eqref{Paradan} are
\begin{align*}
\mp \bfone^w_{\bfC^\sharp_{\widetilde{\Delta}}} &  & \text{and} & & 0.
\end{align*}
Indeed, $\beta(E,\Delta)\cap \partial P=\Delta \cap \partial P$ and
$
\beta(\beta(\varepsilon_i, \widetilde{\Delta}),\Delta) -
\beta(\varepsilon_i,\widetilde{\Delta})=r_i \alpha_{\Delta,k}
$
for $i=1,2$, with $r_1<0$ and $r_2>0$ (cf. Figures~\ref{fi:cubo} and
~\ref{fi:cones3d}).

Consequently, the differences of the contributions of $\Delta$ to the formula in
\eqref{Paradan} before and after $\varepsilon_t$ crosses the wall,
and those of $\widetilde{\Delta}$, sum to zero.

Moreover, for a given $\widetilde{\Delta}\in \PB$, if
$\varphi(\varepsilon_t, \tilde\Delta)$ changes when crossing $E$,
the intersection of $\widetilde\Delta$ with $E$ contains $\Delta \cap P$ for (exactly) one
element $\Delta$ of $\PB$ with $\dim \Delta= \dim \widetilde\Delta -1$ and the
result follows.
\end{proof}

\begin{Remark}\label{re:compare-wpd} This new polytope
decompositions \eqref{Paradan} generalize the weighted version of the Lawrence-Varchenko relation for a simple polytope presented
in \cite{A03}. There, the edge vectors emanating from each vertex are flipped in a systematic way using a polarizing vector, and the weighted characteristic function of the polytope is expressed (\emph{only}) in terms of the weighted characteristic functions of the polarized cones supported at the vertices. In \eqref{Paradan}, not only the polarization is carried
out differently, but, for some values of $\varepsilon$, we consider the weighted characteristic functions of  polarized tangent cones to faces other than vertices. 
Indeed, given $\varepsilon\in W$, we obtain a \emph{different} polarizing vector
for each face of the polytope by taking $\varepsilon$ as starting point, and its projections onto the faces of the polytope as end points, whenever these projections are nonempty. Then  we polarize the tangent cones of the corresponding faces accordingly.  
\end{Remark}

\subsection{Other decomposition formulas}
\begin{figure}[h]
   \centering
  \psfrag{a}{$\varepsilon$}
  \psfrag{b}{$v_0$} 
  \psfrag{c}{$\alpha_{A,1}$} 
  \psfrag{d}{$A$} 
  \psfrag{e}{$B$} 
  \psfrag{f}{$\alpha_{B,1}$} 
  \psfrag{g}{$C$} 
  \psfrag{h}{$\alpha_{C,1}$}
  \psfrag{i}{$H_{v_0,1}^0$} 
  \psfrag{j}{$H_{v_0,2}^0$} 
  \psfrag{k}{$H_{A,1}^0$}
  \psfrag{l}{$H_{B,1}^0$} 
  \psfrag{m}{$H_{C,1}^0$} 
   \includegraphics[scale=.40]{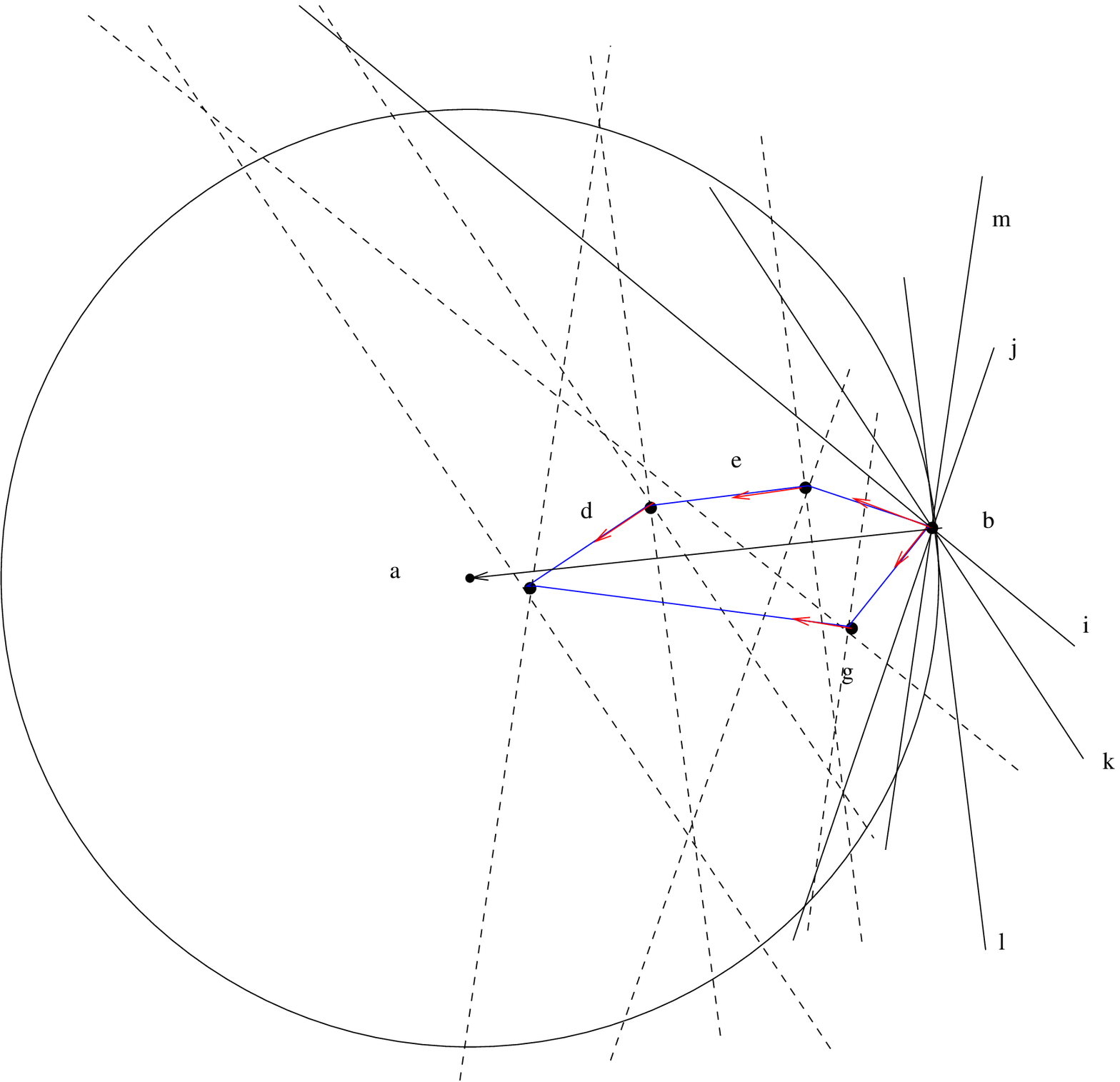}\\
   \caption{}
   \label{fig:figproof}
 \end{figure}

If we polarize the generators of tangent cones with respect to
$\varepsilon-\beta(\varepsilon,\Delta)$ instead of
$\beta(\varepsilon,\Delta)-\varepsilon$, and multiply each term on
the right hand side of \eqref{Paradan} by a factor
$(-1)^{\dim\Delta}$, we obtain new polytope decompositions, under the same hypotheses and
statements of Theorem~\ref{thm:1}:
\begin{Theorem}\label{th:WBG+}
For every compact convex simple polytope $P$ of dimension $d$ in $\R^d$ and
for any $\varepsilon \in W$, we have
\begin{equation}
\label{eq:Paradan-Brianchon-Gram}%
\bfone^w_{P} = \sum_{\Delta \in \PB}(-1)^{m_{\Delta}+\dim\Delta}
\varphi(\varepsilon,\Delta)\bfone^w_{\bfC^\sharp_{\Delta}},
\end{equation}
where the sum is taken over the set $\PB$ of affine spaces generated
by the faces of $P$, where $\bfC^\sharp_\Delta$ is the polarized
tangent cone of $P$ at $\Delta\in \PB$ with respect to the vector
$\varepsilon - \beta(\varepsilon,\Delta)$ (where
$\beta(\varepsilon,\Delta)$ is the orthogonal projection of
$\varepsilon$ onto $\Delta$), where $m_\Delta$ is the number of
generators of the cone $\bfC_\Delta$ whose sign changes by
polarization, and where $\bfone^w_{P}$ and
$\bfone^w_{\bfC^\sharp_{\Delta}}$ are the weighted characteristic
functions of the polytope $P$ and of the polarized cones
respectively.
\end{Theorem}
\begin{proof}
The fact that the right hand side of
\eqref{eq:Paradan-Brianchon-Gram} does not depend on $\varepsilon$
can be proved as in the proof of
Theorem~\ref{thm:1}. Hence, we just have to show that we can find an $\varepsilon$ in some Paradan region 
for which \eqref{eq:Paradan-Brianchon-Gram} holds. 

For that, let us choose an $\varepsilon$ such that $\varphi(\varepsilon, \Delta)=0$ for every $\Delta\in \PB$ with $\dim{\Delta}> 0$, i.e. we choose 
$$
\varepsilon \in \bigcap_{F \,\, \text{a facet of} \,\, P} \left( F + \Delta_F^\perp \right)^c.
$$
Then formula~\eqref{eq:Paradan-Brianchon-Gram} becomes
\begin{equation}\label{eq:Paradan-Brianchon-Gram2}
\bfone^w_{P} = \sum_{v\,\,\text{a vertex of}\,\, P}(-1)^{m_{v}} \bfone^w_{\bfC^\sharp_{v}}.
\end{equation}
Choosing a vector $\xi\in \R^d$ such that, for every vertex $v$ of $P$, $\langle \xi , \alpha_{v,j} \rangle > 0$ whenever $\langle \varepsilon - v , \alpha_{v,j} \rangle > 0$, where the  $\alpha_{v,j}$'s are the edge vectors at $v$, formula~\eqref{eq:Paradan-Brianchon-Gram2} becomes a weighted version of  the Lawrence-Varchenko polytope decomposition (see \cite{L}, \cite{V} \cite{KSW} and \cite{A03}), where the tangent cones at vertices are polarized according to $\xi$, and the result follows\footnote{This weighted version of the  Lawrence-Varchenko relation is different from the ones in \cite{KSW} and \cite{A03} because here we may assign different weights $q_i\in \C$ to the faces of the polytope instead of a fixed complex number. Nevertheless, the proof of this decomposition formula follows easily from the ones in \cite{KSW} and \cite{A03} by applying Lemma~\ref{le:qi's} to the boundary points, as we do in the proof of Theorem~\ref{thm:1}.}.

This choice of polarizing vector $\xi$  can be done in the following way: first we consider the vertex $v_0$ of $P$ that is furthest away from $\varepsilon$; clearly, for $v_0$ we have $\bfC_{v_0}=\bfC_{v_0}^\sharp$ (where this cone is polarized with respect to the vector $\varepsilon - v_0$); then, for any other vertex $v$ and for each edge vector  $\alpha_{v,j}$ satisfying   $\langle \varepsilon - v , \alpha_{v,j} \rangle > 0$, we take the hyperplane  $H_{v,j}^0$ through $v_0$ which is perpendicular to $\alpha_{v,j}$; these hyperplanes intersect at $v_0$ and each of them separates the whole space $\R^d$ into two open regions. Let us denote by  $(H_{v,j}^0)^+$ those regions that contain $\varepsilon$ and take a vector $\xi$ starting at $v_0$ and ending somewhere on the intersection
$$
\bigcap_{v\,\,\text{a vertex of}\,\, P} \bigcap_{\substack{j\,\,\text{s.t.} \\ \langle \varepsilon - v , \alpha_{v,j} \rangle > 0}} (H_{v,j}^0)^+
$$
(we can take for instance $\xi= \varepsilon -v_0$); then clearly  $\langle \xi , \alpha_{v,j} \rangle > 0$ for all edge vectors $\alpha_{v,j}$ satisfying   $\langle \varepsilon - v , \alpha_{v,j} \rangle > 0$ (see Figure~\ref{fig:figproof}). 
\end{proof}

\begin{Remark}\label{rmk:general} We have seen in the above proof that, choosing $\varepsilon$ in an appropriate region, the polytope decomposition formula \eqref{eq:Paradan-Brianchon-Gram} becomes the Lawrence-Varchenko relation. In addition, in some cases, we can also choose $\varepsilon$ so that \eqref{eq:Paradan-Brianchon-Gram} becomes the weighted
Brianchon-Gram formula of \cite{A05}. Indeed, considering for each vertex $v$ of $P$, the cone $\bfC^d_v$ 
generated by the inward normal vectors to the facets through $v$, and taking
the intersection 
$$ P_d:= \bigcap_{v \,
\text{vertex of} \, P} \bfC_v^d,$$
then, whenever $\text{int}
(P_d\cap P)\neq\varnothing$, we can take  $\varepsilon \in
\text{int} (P_d\cap P)$, and obtain $m_\Delta=0$ and
$\varphi(\varepsilon,\Delta)=1$ for every $\Delta$ in $\PB$. 
Then, with this choice of $\varepsilon$, \eqref{eq:Paradan-Brianchon-Gram} becomes the weighted
Brianchon-Gram formula:
\begin{equation}
\label{eq:Brianchon-Gram}%
\bfone^w_{P} = \sum_{F}(-1)^{\dim F}\bfone^w_{\bfC_{F}},
\end{equation}
where the sum is over all faces $F$ of $P$. 
\end{Remark}

\section{The weighted Euler-Maclaurin formula}
\label{se:EML} As an application of our polytope
decompositions, we will give new weighted Euler-Maclaurin formulas
with remainder for the sum of the values of a smooth function $f$ on
the integral points of a simple polytope $P$.
\subsection{Weighted Euler-Maclaurin for intervals}
Let us first recall the weighted Euler-Maclaurin formula for this
sum presented in \cite{AW} (see also \cite{Kn} and \cite{KSW}): let
$q$ be any complex number and  let $f$ be any $\mathcal{C}^m$
function on the real line ($m\geq 1$); for integers $a<b$, the sum
\begin{eqnarray}\label{eq:emlinterval}
\sum_{[a,b]}{}^{q} \, f &  := &   q f(a) + f(a+1) + \cdots + f(b-1) +  q f(b) \\
 \nonumber & = & \bfQ^{2k}_q(D_1)\, \bfQ^{2k}_{q}(D_2) \int_{a-h_1}^{b+h_2} f(x)\,
 dx\, \bvert_{h_1=h_2=0} \, + R_{m}(f),
\end{eqnarray}
where $k= \lfloor m/2 \rfloor$, where $\bfQ^{2k}_q(S)$ denotes the
truncation  at the even integer $2k$ of the power series
\begin{equation}
\bfQ_q(S) = (q-1)S + \text{{\bf Td}}(S) = 
 1 + \left(q-\frac{1}{2}\right)S+\sum_{k=1}^\infty
\frac{b_{2k}}{(2k)!} \, S^{2k} = \left(q-\frac{1}{2}\right)S +
\frac{S/2}{\mathrm{tanh}(S/2)}
\end{equation}
(here $\text{{\bf Td}}$ is the classical Todd function defined by
$\text{{\bf Td}}(S):=S/(1-e^{-S})=1-b_1S+\sum_{k=1}^\infty
\frac{b_{2k}}{(2k)!} \,  S^{2k}$, with  $b_k$ the $k$-th Bernoulli
number \cite{B}), where
\begin{equation*}
D_1:=\frac{\partial}{\partial h_1}, \quad  D_2:=\frac{\partial}{\partial h_2},
\end{equation*}
and
\begin{equation}
R_{m}(f) : =(-1)^{m-1} \int_a^b P_{m}(x)f^{(m)}(x)\, dx,
\end{equation}
with
\begin{equation}
P_{2k+1} :=  (-1)^{k-1} \sum_{n=1}^\infty\frac{2 \sin(2n\pi
x)}{(2n\pi)^{2k+1}} = \frac{1}{(2k+1)!} B_{2k+1}(\{x\})
\end{equation}
if $m=2k+1$ is odd,
and
\begin{equation}
P_{2k} :=  (-1)^{k-1} \sum_{n=1}^\infty\frac{2 \cos(2n\pi
x)}{(2n\pi)^{2k}} = \frac{1}{(2k)!} B_{2k}(\{x\})
\end{equation}
if $m=2k$ is even, (here $B_{m}$ is the $m$-th Bernoulli polynomial
and $\{x\}:=x-\lfloor x\rfloor$ is the fractional part of $x$).
\begin{Remark}
The functions $\bfQ^{2k}_q$ satisfy the following symmetry property
\begin{equation}\label{eq:symQ}
\bfQ^{2k}_q(S)=\bfQ^{2k}_{1-q}(-S).
\end{equation}
Indeed, $\bfQ^{2k}_q(S)$ is a polynomial with constant coefficients,
$1+(q-\frac1{2})S +$ terms of even degree independent of $q$.
\end{Remark}
Equation \eqref{eq:emlinterval}, when applied to a $\mathcal{C}^{m}$
function of compact support, gives the weighted Euler-Maclaurin
formula for the half ray $[a,\infty)$:
\begin{equation}\label{eq:emlhalfray+}
\sum_{[a,\infty)}{}^{q} \, f  : =   q f(a) + f(a+1) + f(a+2) + \cdots =  \bfQ^{2k}_q(D_1) \int_{a - h_1}^{\infty} f(x)\, dx \bvert_{h_1=0} \, + R_{m}(f),
\end{equation}
where
\begin{equation}
R_{m} : =(-1)^{m-1} \int_a^\infty P_{m}(x)f^{(m)}(x)\, dx.
\end{equation}
Moreover, for the  half ray $(-\infty, a]$, we have
\begin{eqnarray}
\sum_{(-\infty,a]}{}^{q} \, f &  = &   q f(a) + f(a-1) + f(a-2) + \cdots
\end{eqnarray}
and so, considering the function $g$ defined by $g(x)=f(a-x)$, we obtain
\begin{eqnarray}\label{eq:emlhalfray-}
 \,\,\,\,\,\,\,\,\,\, \sum_{(-\infty,a]}{}\!\!\!\!^{q}\! f \!\!\!\!\! & =  & \!\!\!\!\!\!\!\!\sum_{[0,\,\infty)}{}\!\!^{q} \, g =  \bfQ^{2k}_q(D_1) \int_{- h_1}^{\infty} g(x)\, dx\, \bvert_{h_1=0} \, +
 (-1)^{m-1} \int_0^\infty P_{m}(x)g^{(m)}(x)\, dx \\
\nonumber & = & \bfQ^{2k}_q(D_1) \int^{a + h_1}_{-\infty} f(x)\, dx
\, \bvert_{h_1=0} \, + R_{m}(f),
\end{eqnarray}
where
\begin{equation}
R_{m}(f) : = (-1)^{m-1} \int^a_{-\infty} P_{m}(x)f^{(m)}(x)\, dx
\end{equation}
(here we used the parity and the  $2\pi$-periodicity of $\sin(x)$
and $\cos(x)$). From formulas  \eqref{eq:emlhalfray+} and
\eqref{eq:emlhalfray-} and symmetry property \eqref{eq:symQ}, we
obtain the Euler-Maclaurin formula for the whole real line $\R$:
\begin{equation}
\label{eq:emlreals} \sum_{\R}{}^{'} f : = \sum_{x\in \Z}
f(x) = \sum_{(-\infty,0]}{}\!\!\!^{q} f  +
\sum_{[0,\infty)}{}\!\!^{(1-q)}  f = \int_\R
f(x) dx + (-1)^{m-1} \! \int_\R P_{m}(x)f^{(m)}(x)\, dx.
\end{equation}
\subsection{Twisted weighted Euler-Maclaurin formulas for intervals}
We will now consider the twisted weighted sum for a half ray
\begin{equation}
\sum_{n \geq 0}{}^q \lambda^n f(n) = qf(0) + \sum_{n = 1}^\infty \lambda^n f(n)
\end{equation}
where $\lambda\neq 1$ is a $K$-th root of unity with $K$ a positive
integer. Let $Q_{m,\lambda}$ be the distributions defined recursively in
\cite{KSW} by
$$Q_{0,\lambda}(x): =-\sum_{n\in \Z}  \lambda^n \delta(x-n)$$
and
$$
\frac{d}{d\, x} Q_{m,\lambda}(x)= Q_{m-1,\lambda}(x) \quad
\text{and} \quad \int_o^K Q_{m,\lambda} (x)\, dx = 0.
$$
Moreover, let us consider the polynomials defined in \cite{AW} by
$$
\bfN_q^{k,\lambda}(S): =\left( q+\frac{\lambda}{1-\lambda} \right) S
+ Q_{2,\lambda}(0)S^2 +  Q_{3,\lambda}(0)S^3+ \cdots +
Q_{k,\lambda}(0)S^k,
$$
where $\lambda\neq 1$ is a root of unity. 

Since $\frac{d}{d x} {\bf 1}_{[n,n+1)} (x)= \delta(x-n) - \delta(x-(n+1))$, we have 
$$
\frac{d}{d x}\left(\sum_{n\in \Z} \lambda^n {\bf 1}_{[n, n+1)} (x) \right) = \frac{\lambda -1}{\lambda} \sum_{n\in \Z} \lambda^n \delta(x-n) =  \frac{1- \lambda}{\lambda} \, Q_{0,\lambda}(x),
$$
implying that $Q_{1,\lambda}(x)= \frac{\lambda}{1- \lambda} \sum_{n\in \Z} \lambda^n {\bf 1}_{[n,n+1)}$ (note that $\int_0^K Q_{1, \lambda} (x) \, dx =  \frac{\lambda}{1 - \lambda} \displaystyle{\sum_{n=0}^{K-1} \lambda^n} =0$). On the other hand, integrating by parts, we have
$$
\int_0^\infty Q_{1, \lambda}(x) f'(x) \, dx =  \frac{\lambda}{1- \lambda} \sum_{n=0}^\infty \int_n^{n+1} \lambda^n\, f'(x)\, dx = - \frac{\lambda}{1- \lambda} f(0) + \lambda f(1) + \lambda^2 f(2) + \cdots 
$$
and so,
\begin{align*}
&q\,f(0) + \sum_{n\geq 1} \lambda^n f(n)\,  =  (q+ \frac{\lambda}{1- \lambda})f(0) + \int_0^\infty Q_{1,\lambda}(x) f'(x)\, dx \\
 & =  (q+ \frac{\lambda}{1- \lambda})f(0) - Q_{2,\lambda}(0) f'(0) +  Q_{3,\lambda}(0) f''(0) - \cdots + (-1)^{k-1} Q_{k,\lambda}(0) f^{(k-1)}(0) + \\ & \,\,\,\,\,\,\,\,\,\,\,\, + (-1)^{k-1} \int_0^\infty  Q_{k,\lambda}(x) f^{(k)}(x) \, dx.
\end{align*}
Then, since $(-1)^m f^{(m-1)} (0)=\left(\frac{\partial}{\partial\, h}\right)^m \int_{-h}^\infty f(x)\, dx_{\vert_{h=0}}$, 
we obtain the following twisted Euler-Maclaurin formula:
\begin{Proposition}(\cite{AW}, \cite{KSW})
Let $k>1$ and let $f\in \mathcal{C}^k(\R)$ be compactly supported. Then
\begin{equation}\label{twistedeml+}
\sum_{n\geq 0}{}^q \lambda^n f(n) = 
\bfN_q^{k,\lambda} \left(\frac{\partial}{\partial h}\right)
\left.\int_{-h}^{\infty} f(x)\, dx \, \right\vert_{h=0} +
(-1)^{k-1}\int_0^\infty Q_{k,\lambda}(x)f^{(k)}(x)\, dx.
\end{equation}
\end{Proposition}

\begin{Remark}\label{rmk:Vergne}
If, for $\lambda \neq 1$, we write $\lambda = e^{2\pi i j/K}$, then, by the Poisson formula, we have 
\begin{align*}
Q_{0,\lambda}(x) & = -\sum_{n \in \Z}\lambda^n\, \delta(x-n) = - \sum_{n\in \Z} \lambda^x \delta(x-n)\\ & = - e^{2\pi i x \frac{j}{K}} \sum_{n\in \Z} \delta(x-n) =  - e^{2\pi i x\frac{j}{K}} \sum_{r\in \Z} e^{2\pi i r\, x} =- \sum_{r\in \Z} e^{2\pi i (r+\frac{j}{K})\,x}. 
\end{align*}
Hence, for $m>1$, we obtain
\begin{equation}\label{eq:sym0}
Q_{m, \lambda}(x)=- \frac{1}{(2\pi i)^m }\sum_{r\in \Z} \frac{ e^{2\pi i (r+\frac{j}{K})\,x}}{(r+\frac{j}{K})^m},
\end{equation}
and so $Q_{m,\lambda}(0)=- \frac{1}{(2\pi i)^m } \sum_{r\in \Z} \frac{1}{(r+\frac{j}{K})^m}$ is the $(m-1)$-th coefficient of the Taylor series expansion of $\frac{1}{1-e^{2\pi i \frac{j}{K}-s}}=\frac{1}{1-\lambda e^{-s}}$ at $s=0$ (the derivative of $\frac{1}{1-e^{2\pi i \frac{j}{K}-s}}$ with respect to $s$ is equal to $\frac{1}{4 \sin^2{(\frac{\pi j}{K} - \frac{s}{2i})}}= \frac{1}{4\pi^2} \sum_{r\in \Z} \frac{1}{(r+\frac{j}{K}-\frac{s}{2\pi i})^2}$,\footnote{Note that $\frac{\pi^2}{\sin^2{\pi z}} =\sum_{r\in \Z} \frac{1}{(r+z)^2}$.} and higher order derivatives are obtained differentiating this series expansion). Consequently, considering the operators
$$
{\bf T}(\lambda,S):= \frac{S}{1-\lambda e^{-S}}
$$
defined in \cite{BrV}, we have that $\bfN_q^{k,\lambda}(S)$ is the truncation at the integer $k$ of the power series
$$
\bfN_q^{\lambda}(S):= (q+\frac{\lambda}{1-\lambda})S - \frac{S}{1-\lambda}+ {\bf T}(\lambda,S) = (q-1) S + {\bf T}(\lambda,S).
$$
\end{Remark}
From \eqref{eq:sym0} it is clear that the operators $\bfN_q^{m,\lambda}$ satisfy
the following symmetry property
\begin{equation}\label{eq:symbfN}
\bfN_{1-q}^{m,\lambda^{-1}} (S)= \bfN_q^{m,\lambda} (-S).
\end{equation}
\begin{Remark}\label{re:multipleofS}
If, for $\lambda=1$, we define
$$
\bfN_q^{k,1}(S): = \bfQ_q^{2\lfloor k/2 \rfloor}(S) \quad \text{and}\quad Q_{k,1}: =P_k,
$$
then formula \eqref{twistedeml+} becomes formula \eqref{eq:emlhalfray+}
and so it is still valid. Note that, if $\lambda \neq 1$,
$\bfN_q^{k,\lambda} (S)$ is a multiple of $S$ and that, if $\lambda =1$,
then $\bfN_q^{k,\lambda}(S)=1+$ a multiple of $S$.  Moreover, still when $\lambda=1$,
symmetry property \eqref{eq:symbfN} becomes property \eqref{eq:symQ}.
\end{Remark}
\subsection{Weighted Euler-Maclaurin formulas for cones}
For a subset $J\!\subset \!\{1,\ldots, d\}$, let $\bfS_J$ be the
standard $J$-sector $\bfS_J:=\{x\in \R^d\mid x_j \geq
0\,\,\text{for}\,\, j\in J\}$. Iterating equations
\eqref{eq:emlhalfray+} and \eqref{eq:emlreals}, we obtain an
Euler-Maclaurin formula for $\bfS_J$ ($J \neq \varnothing$) and a $\mathcal{C}^{m}$
function of compact support:
\begin{equation}\label{eq:emlsector}
\sum_{\bfS_J \, \cap \, \Z^d}{}\!\!^{w} \, f   :=\!\!\!\!\!\!\!\!\!
\sum_{\substack{ x_j \in \Z^+\!\!, j \in J \\[.5mm]  x_j \in \Z,\, j \notin J}} (\bfone_{\bfS_J}^w  f)
(x_1, \ldots, x_d) = \prod_{j\in J} \bfQ^{2k}_{q_j}(D_j) \int_{\bfS_J(h_J)}  f(x) \, dx \bvert_{h_J = 0}
\,  + R_{m}^{J_\mathrm{st}}(f),
\end{equation}
where $\bfone^w_{\bfS_J}$ is the weighted characteristic function
for the $J$-sector defined in Remark~\ref{le:wcharcone},  where
$D_i=\partial/ \partial h_i$, where $h_J=(h_{j_1},\ldots,h_{j_n})$
with $J=\{j_1,\ldots, j_n\}$, where $\bfS_J (h_J)=\{x\in \R^d\mid
x_j \geq - h_j,\,\,\text{for}\,\,j\in J\}$ is the shifted
$J$-sector, and where the remainder  $R_{m}^{J_\mathrm{st}}(f)$ is
given by
\begin{align} \nonumber
R_{m}^{J_\mathrm{st}}(f) & : =\\ \nonumber = \sum_{\substack{I  \subseteq J}}\!\!\!\!\!\!
\sum_{\substack{R \supseteq J\\ R\subseteq \{ 1,\ldots,d\} \\ R \neq
I }}\!\!\!\!\!\! & (-1)^{(m-1)(\vert R \vert-\vert I\vert)}
\prod_{i\in I} \bfQ^{2k}_{q_i}(D_i) \int_{\bfS_J(h_J)} \prod_{i
\in R\setminus I} \!\!\!\! P_{m}(x_i)\!\!\!\!  \prod_{j \in R\setminus I}\!\!\!\! \left(
\frac{\partial}{\partial x_j} \right)^{\!\!m}\!\!\!\!  f(x) \, dx
\bvert_{h_J=0}.
\end{align}
If $J= \varnothing$ then $\bfS_J$ is the whole space $\R^d$ and so
\begin{eqnarray}\label{eq:emlsector2}
\sum_{\bfS_J \, \cap \, \Z^d}{}\!\!^{w} \, f  = \int_{\R^d}  f(x) \,
dx + R_{m}^{\varnothing}(f),
\end{eqnarray}
with
\begin{eqnarray}
\quad\quad\quad R_{m}^{\varnothing}(f) : = \sum_{\substack{R \neq
\varnothing \\ R \subseteq \{ 1,\ldots,d\} \\ R \neq I }} (-1)^{(m-1)
\vert R \vert} \int_{\R^d}  \prod_{i \in R} P_{m}(x_i) \prod_{j \in
R} \left( \frac{\partial}{\partial x_j} \right)^{\!\!m}  f(x) \, dx.
\end{eqnarray}
Let us now consider a regular integral $J$-sector $\bfC_J$, the image of the standard
$J$-sector by an affine transformation
$$
x \mapsto A_{_J}x:=Mx + b, \quad \text{with} \quad M\in SL(d,\Z) \quad \text{and} \quad b\in \R^d.
$$
Moreover, let us denote by $\bfC_J(h)$ the expanded sector, image of
$\bfS_J(h)$ under this affine transformation. For a
$\mathcal{C}^{m}$ function of compact support $f$, let us consider
$g:= A_J^*f =f \circ A_J$. Then,
\begin{eqnarray*}
\sum_{\bfC_J \, \cap \, \Z^d}{}\!\!^{w} \, f := \sum_{\bfS_J \, \cap
\, \Z^d}{}\!\!^{w} \, g = \prod_{j\in J} \bfQ^{2k}_{q_j}(D_j)
\int_{\bfS_J(h_J)}  g(x) \, dx \bvert_{h_J = 0} \,  +
R_{m}^{J_\mathrm{st}}(g),
\end{eqnarray*}
and we obtain the following Euler-Maclaurin formula for a regular
$J$-sector:
\begin{equation}
\label{eq:emlJsector} \sum_{\bfC_J \, \cap \, \Z^d}{}\!\!^{w} \, f =
\prod_{j\in J} \bfQ_{q_j}^{2k}(D_j) \int_{\bfC_J(h_J)}  f(x) \, dx
\bvert_{h_J = 0} \,  + R_{m}^{\bfC_J}(f),
\end{equation}
where $R_{m}^{\bfC_J}(f)= R_{m}^{J_\mathrm{st}}(g)$. 

\subsection{Weighted Euler-Maclaurin formula for regular simple integral polytopes}
From \eqref{eq:emlJsector} we can write an  Euler-Maclaurin formula for a
regular integral polytope $P$ with $N$ facets, by using a
polytope decomposition from Theorem~\ref{thm:1}:
\begin{align}
\label{eq:emlregpolytope}
\begin{split}
& \sum_{P \, \cap \, \Z^d}{}\!\!^{w} \, f := \sum_{P \, \cap \,
\Z^d} \bfone_P^w f = \sum_{\Delta \in \PB} (-1)^{m_\Delta}
\varphi(\varepsilon,\Delta) \sum_{\bfC^\sharp_{\Delta} \, \cap \, \Z^d}{}\!\!^{w} \, f \\
& = \sum_{\Delta \in \PB} (-1)^{m_\Delta}
\varphi(\varepsilon,\Delta) \left( \prod_{j=1}^N \bfQ^{2k}_{q_j}(D_j)
\int_{\bfC^\sharp_{\Delta}(h_\Delta)}  f(x) \, dx \bvert_{h_\Delta = 0}  +
R_{m}^{\bfC^\sharp_{\Delta}}(f)   \right) \\
& = \prod_{j=1}^{N} \bfQ_{q_j}^{2k}(D_j) \int_{P(h_1,\ldots, h_N)}
f(x) \, dx \bvert_{h = 0} \!\!  + S_{m}^{P}(f),
\end{split}
\end{align}
where
\begin{equation}\label{eq:finalremainder}
S_{m}^{P}(f):= \sum_{\Delta \in \PB}  (-1)^{m_\Delta}
\varphi(\varepsilon,\Delta)  R_{m}^{\bfC^\sharp_{\Delta}}(f),
\end{equation}
and where the dilated polytope $P(h_1,\ldots, h_N)$ is obtained by
shifting the $i$th facet outward by a ``distance'' $h_i$. Here we
used the fact that, when multiplying the differential operator in the
first term of the right hand side of \eqref{eq:emlJsector} by any
operator of the form $\bfQ^{2k}_{q_j}(D_j)$, with $j\notin J$, all
that will remain of $\bfQ^{2k}_{q_j}(D_j)$ is the constant term $1$,
not affecting the final result. Note that both $\sum_{P \, \cap \,
\Z^d}{}\!\!^{w} \, f$ and  $\prod_{j=1}^N \bfQ_{q_j}^{2k}(D_i)
\int_{P(h_1,\ldots,h_N)}  f(x)$ do not depend on the choice of
$\varepsilon$ (that is, do not depend on the Paradan region used). Consequently, the remainder is also independent
of this choice.
\begin{Remark} Alternatively, using the polytope decomposition of Theorem~\ref{th:WBG+} we obtain
a different expression for the remainder in
\eqref{eq:finalremainder}. Indeed, we get
\begin{equation}\label{eq:finalremainder2}
S_{m}^{P}(f):= \sum_{\Delta \in \PB}  (-1)^{m_\Delta+\dim{\Delta}}
\varphi(\varepsilon,\Delta)  R_{m}^{\bfC^\sharp_{\Delta}}(f),
\end{equation}
where now the tangent cones $\bfC_\Delta$ are polarized with respect to the
vectors $\varepsilon-\beta(\varepsilon,\Delta)$.
\end{Remark}
\subsection{Weighted Euler-Maclaurin formula for simple integral polytopes}
To extend formula \eqref{eq:emlregpolytope} to simple integral
polytopes we need to obtain an Euler-Maclaurin formula for simple
$J$-sectors. We can describe a simple   $J$-sector, $\bfC_J$,  with
$J\!=\{i_1,\ldots, i_n\}\!\subset \{1,\ldots, d\}$ as the intersection
of $n$ half spaces $H_i$  in general position
\begin{equation}\label{eq:C_J}
\bfC_J= \bigcap_{j \in J} H_j,
\end{equation}
where $H_j:= \{x\in \R^d\mid \langle x, \eta_j  \rangle + \lambda_j
\geq 0 \}$ for rational vectors $\eta_j$. Clearly, $\bfC_J$ is a
cone of the form $\bfC_{\Delta_J}^\sharp$ along the affine space $\Delta_J \in \PB$
defined by
\begin{equation*}
\Delta_J=  \bigcap_{j \in J}\, \partial H_j
\end{equation*}
(see \eqref{eq:conealongDelta}). Clearing denominators we can assume
the $\eta_i$'s to be integral and we impose the normalizing condition
that they are primitive elements of the dual lattice $\Z^n{}^*$
(note that these vectors are inward normals to the facets of
$\bfC_J$).  Let us take the dual basis $\{\alpha_{i_1},\ldots,
\alpha_{i_n}\}$ in $\R^n$ (that is, such that $\langle \alpha_{k} ,
\eta_{l} \rangle=\delta_{k l}$ for $k,l\in J$) and denote by $T_J
\subseteq \R^{d}{}^*$ the  subspace generated by these vectors. The
$\alpha_i$'s are what, in Section~\ref{se:decompositions}, we called
generators of $\bfC_{\Delta_J}$  and generate a lattice $\hat{\ell}$
in $\R^n{}^*$ which is a finite extension of $\Z^n$ (this extension
is trivial exactly when $\bfC_J$ is regular). Let $\Gamma_J$ be the
finite group
$$
\Gamma_J= (\Z^n{}^* \cap T_J) / \hat{\ell}^*.
$$
This group  is trivial exactly when $\bfC_J$ is regular, and its
order is $ \vert\Gamma_J \vert = \vert \det{\hat{\ell}} \vert$.

Moreover, as it is shown in \cite{KSW}, $\gamma \mapsto e^{2\pi i
\langle \gamma, x \rangle}$ defines a character of $\Gamma_J$,
whenever $x\in \hat{\ell}$, which is trivial iff $x\in \Z^n$. Since,
by a theorem of Frobenius, the average value of a character on a
finite group is equal to zero if the character is non-trivial, and equal to one
otherwise, we have
\begin{equation*}
\frac{1}{\vert\Gamma_J \vert}\sum_{\gamma \in \Gamma_J}  e^{2\pi i
\langle \gamma, x \rangle}= \left\{ \begin{array}{ll} 1 & \text{if}
\quad x \in \Z^n  \\ 0  & \text{otherwise}
\end{array} \right.,
\end{equation*}
for all $x \in \hat{\ell}$. Consequently, for any compactly
supported function $f$ on $\R^n$,
\begin{equation}\label{eq:emlsimple1}
\sum_{\bfC_J  \, \cap \, \Z^d}{}\!\!\!^{w} \, f =
\frac{1}{\vert\Gamma_J \vert} \sum_{\gamma \in \Gamma}
\sum_{x}\!^{w}  e^{2\pi i \langle \gamma,x\rangle} f(x),
\end{equation}
where we sum over all
\begin{equation}\label{eq:x}
x= y + \sum_{j\in J} m_j \alpha_{j}
\end{equation}
with $y \in \Delta_J \cap \Z^n$ and all $m_j\in\Z^+$ with $j\in J$.
Moreover, the cone $\bfC_J$ is the image of the standard $J$-sector
$\bfS_J$ under an affine map
\begin{equation}\label{eq:changevar}
t \mapsto A_{_J}t:= U_J t + b, \quad \text{with} \quad b\in \R^d,
\end{equation}
where, for $j\in J$,  $U_J\in GL(d,\Z)$ carries the vectors $e_j$ of the
standard basis of $\R^d$ into the basis $\{\alpha_{i_k}\}_{k=1}^n$.
Hence, $\vert \det U_J \vert = 1/\vert \det{\hat{\ell}}\vert = 1/\vert
\Gamma_J \vert$. On the other hand, since in \eqref{eq:x} we have  $y\in
\Z^n$, we get
$$
e^{2\pi i\langle \gamma, x \rangle} = \prod_{j\in J}
\lambda_j^{m_j},\quad \text{with} \quad \lambda_j= e^{2\pi i\langle
\gamma, \alpha_j \rangle},
$$
and so the inner sum in \eqref{eq:emlsimple1} becomes
\begin{eqnarray}\nonumber
\lefteqn{ \sum_{j\in J} \, \sum_{y\in \Delta_J \cap \Z^n} \,  \sum_{m_{j} \geq 0}{}\!\!^{q_{j}} \, \left( \prod_{l\in J} \lambda_{l}^{m_{l}} \right) \,
f(y+\sum_{l\in J}  m_{l}\alpha_{l}) =} \\ \label{eq:emlsimple2} 
& = &  \sum_{j\in J} \, \sum_{\substack{1 \leq  i \leq d \\[0.5 mm] i \notin J}}
\, \sum_{m_i \in \Z} \,  \sum_{m_{j} \geq 0}{}\!\!^{q_{j}} \, \left(\prod_{l\in J} \lambda_{l}^{m_{l}} \right) \, g(m_{1},\ldots, m_{d}),
\end{eqnarray}
where $g=f \circ A_J $. Iterating the twisted remainder formula for
the half ray \eqref{twistedeml+} and the Euler-Maclaurin formula
\eqref{eq:emlreals} for the whole real line,  the sum in
\eqref{eq:emlsimple2} can be written as
\begin{equation}\label{eq:emlsimple3}
\prod_{j\in J} \,
N_{q_j}^{k,\lambda_j}\left(\frac{\partial}{\partial h_j}\right)
\int_{\bfS_J(h_{j_1},\ldots,h_{j_n})}  g_J(t) \, dt\bvert_{h=0}\,
+\, R_{{\bf
q}_J,k}^{\text{std}}(\lambda_{j_1},\ldots,\lambda_{j_n};g),
\end{equation}
where again $\bfS_J(h_{j_1},\ldots,h_{j_n})=\{ (t_1,\ldots,t_d)\mid
t_i\geq - h_i\quad \text{for}\quad j\in J\}$ denotes the dilated
standard $J$-sector, and where, for ${\bf q}_J:=(q_{j_1}, \ldots, q_{j_n})$, the remainder  is given by
\begin{eqnarray}\label{eq:remainder3}
& & \lefteqn{\quad R_{{\bf q}_J,k}^{\text{std}}(\lambda_{\gamma,
j_1},\ldots,\lambda_{\gamma, j_n};g):= \sum_{\substack{I \subset J}}
\sum_{\substack{R\supseteq J\\ R\subseteq \{1,\ldots, d \} \\ R \neq
I }}(-1)^{(k-1)(\vert R \vert - \vert I \vert)} } \\ \nonumber & &
\prod_{i\in I} N_{q_i}^{k,\lambda_{\gamma,
i}}\left(\frac{\partial}{\partial h_i}\right)
\int_{\bfS_J(h_{j_1},\ldots,h_{j_n})} \prod_{j\in R \setminus I}
Q_{k,\lambda_{\gamma, j}}(t_j) \prod_{j\in R \setminus I}
\left(\frac{\partial}{\partial t_j}\right)^k \, g(t) \,  dt \,
\bvert_{h=0},
\end{eqnarray}
with $g = f\circ A$. Changing variables by the inverse
transformation of \eqref{eq:changevar}, the Euler-Maclaurin formula
in \eqref{eq:emlsimple1} becomes
\begin{equation}\label{eq:emlsimple4}
\sum_{\bfC^J  \, \cap \, \Z^d}{}\!\!\!^{w} \, f = \sum_{\gamma \in
\Gamma} \,\, \prod_{j\in J} \,
N_{q_j}^{k,\lambda_{\gamma,j}}\left(\frac{\partial}{\partial
h_j}\right) \int_{\bfC_J(h_J)}  f(x)\,  dx \bvert_{h_J=0} + R_{{\bf
q}_J,k}^{\bfC_J}(f),\end{equation} where $\lambda_{\gamma,j}:=e^{2
\pi i \langle \gamma, \alpha_j\rangle}$, where, for
$h_J:=(h_{j_1},\ldots,h_{j_n})$, $\bfC_J(h_J)$ denotes the image of
the dilated standard $J$-sector $\bfS_J(h_{J})$ under the affine
transformation $A_J$ defined in  \eqref{eq:changevar}, and where the
remainder $R_{{\bf q}_J,k}^{\bfC_J}(f)$ is given by
\begin{align}\label{eq:remainder4}
\begin{split}
R_{{\bf q}_J, k}^{\bfC_J}(f):= & \sum_{\gamma \in \Gamma}
\sum_{\substack{I \subset J}}   \sum_{\substack{R\supseteq J\\
R\subseteq \{1,\ldots, d \} \\ R \neq I }}(-1)^{(k-1)(\vert R \vert
-\vert I\vert)}  \prod_{i\in I}
N_{q_i}^{k,\lambda_{\gamma,i}}\left(\frac{\partial}{\partial h_i}\right)\\
& \int_{\bfC_J(h_J)}  \prod_{j\in R \setminus I}
Q_{k,\lambda_{\gamma,j}}((U^{-1}_{jk})_{k} (x-b) ) \prod_{j\in R
\setminus I} D_{j}^k f(x) \,  dx \bvert_{h=0},
\end{split}
\end{align}
with $(U^{-1}_{jk})_{k}$ the $j$th row of $U^{-1}$ and with $D_{j}$
the directional derivative along the $j$th column vector of
$U^{-1}$. Note that, when $j \in J$, this is the directional
derivative along $\alpha_j$.

Let now $P$ be  a simple polytope and choose an $\varepsilon$ on
some Paradan region.  For each affine space $\Delta$ generated by a
face of $P$ there is a $J$-sector $\bfC_J$ equal to the polarized
tangent cone $\bfC_\Delta^\sharp$ of $\Delta$
(cf. \eqref{eq:polarizedcone}) and so we can associate a finite group
$\Gamma_\Delta$ to $\Delta$ by simply taking the corresponding group
$\Gamma_J$.   Let $P(h)$ denote the dilated polytope obtained by
shifting the $i$th facet by a distance $h_i$. Our decompositions
of $P(h)$ involve dilated sectors but now, dilating
the facets of $P$ outward results in dilating some of the facets of
$\bfC_\Delta^\sharp$ inward and some outward. Explicitly, taking
$J\subset \{1,\ldots,d\}$ such that $\bfC_\Delta^\sharp=\bfC_J$ (see
\eqref{eq:C_J}), the inward normal vector to the $j$th facet of
$\bfC_\Delta^\sharp$  ($j\in J$) is
\begin{equation*}
\eta_{\Delta,j}^\sharp = \left\{ \begin{array}{ll} \eta_j, &
\text{if} \quad \alpha_{\Delta,j}^\sharp=\alpha_{\Delta,j} \\
-\eta_j, & \text{if} \quad
\alpha_{\Delta,j}^\sharp=-\alpha_{\Delta,j}. \end{array}\right. ,
\end{equation*}
where $\eta_j$ is the inward pointing primitive normal vector to the
$j$th facet of $P$ (note that $\alpha_{\Delta,j}^\sharp$, $j\in J$,
is  the dual basis to the corresponding vectors $\eta_j$). The
dilated sectors that appear on the right side of the polytope decompositions of  $P(h)$ are then
$\bfC_\Delta^\sharp(h_{\Delta,j_1}^\sharp,\ldots,
h_{\Delta,j_n}^\sharp)$ (with $J=\{j_1,\ldots,j_n\}$), where
\begin{equation*}
h_{\Delta,j_i}^\sharp = \left\{ \begin{array}{ll} h_{j_i}, &
\text{if} \quad \alpha_{\Delta,j_i}^\sharp=\alpha_{\Delta,j_i} \\
-h_{j_i}, & \text{if} \quad
\alpha_{\Delta,j_i}^\sharp=-\alpha_{\Delta,j_i}. \end{array}\right.
\end{equation*}
Moreover, the roots of unity that appear in the Euler-Maclaurin
formula for $\bfC_\Delta^\sharp$ are
\begin{equation*}
\lambda_{\gamma,j,\Delta}^\sharp = e^{2 \pi i \langle \gamma,
\alpha_{\Delta,j}^\sharp \rangle} = \left\{ \begin{array}{ll}
\lambda_{\gamma,j,\Delta}, & \text{if} \quad
\alpha_{\Delta,j}^\sharp=\alpha_{\Delta,j} \\
\lambda_{\gamma,j,\Delta}^{-1}, & \text{if} \quad
\alpha_{\Delta,j}^\sharp=-\alpha_{\Delta,j}. \end{array}\right.
\end{equation*}

Hence, for any compactly supported function in $\R^d$ of type
$\mathcal{C}^{nk}$ (for an integer $k\geq 1$) the 
decomposition formula of Theorem~\ref{thm:1} applied to $P(h)$  along with  formula \eqref{eq:emlsimple4} give
\begin{align}\label{eq:emlsimplepoly}
\sum_{P \cap \Z^d}{}^{w} \, f &= \sum_{\Delta \in \PB}
(-1)^{m_\Delta} \varphi(\varepsilon,\Delta)
\sum_{\bfC^\sharp_{\Delta}\cap \Z^n}{}\!^{w} \, f = \\ \nonumber
\sum_{\Delta \in \PB} (-1)^{m_\Delta}
\varphi(\varepsilon,\Delta)& \sum_{\gamma \in \Gamma_\Delta} \!\!\! \!\!\! \prod_{\substack{j\in J_\Delta \\[.5mm]
J_\Delta = \{ j_1, \ldots, j_{n_\Delta}\} }} \!\!\! \!\!\!\!\!\!{\bf
N}^{k,\lambda^\sharp_{\gamma,j,\Delta}}_{q^\sharp_{j}}
\left(\frac{\partial}{\partial h^\sharp_{\Delta,j}}\right)
\int_{\bfC^\sharp_{\Delta}(h^\sharp_J)} f(x) \, dx
\bvert_{h_J^\sharp=0}  \!\!\! + R^P_{w,k}(f),
\end{align}
where $h_J^\sharp = (h^\sharp_{i_1},\ldots,
h^\sharp_{i_{n_\Delta}})$ and where
\begin{equation}\label{eq:remaindersimplepoly} R^P_{w,k}(f) :=
\sum_{\Delta \in \PB} (-1)^{m_\Delta} \varphi(\varepsilon,\Delta)
R^{\bfC^\sharp_{\Delta}}_{q_{J_\Delta}^\sharp, k}(f).
\end{equation}

\begin{Remark}
Using the polytope decomposition of Theorem~\ref{th:WBG+} we obtain
\begin{align}\label{eq:emlsimplepolyB}
\sum_{P \cap \Z^d}{}^{w} \, f &= \sum_{\Delta \in \PB}
(-1)^{m_\Delta+\dim{\Delta}} \varphi(\varepsilon,\Delta)
\sum_{\bfC^\sharp_{\Delta}\cap \Z^n}{}\!^{w} \, f = \\ \nonumber
\sum_{\Delta \in \PB} (-1)^{m_\Delta+ \dim{\Delta}}
\varphi(\varepsilon,\Delta)& \sum_{\gamma \in \Gamma_\Delta} \!\!\! \!\!\! \prod_{\substack{j\in J_\Delta \\[.5mm]
J_\Delta = \{ j_1, \ldots, j_{n_\Delta}\} }} \!\!\! \!\!\!\!\!\!{\bf
N}^{k,\lambda^\sharp_{\gamma,j,\Delta}}_{q^\sharp_{j}}
\left(\frac{\partial}{\partial h^\sharp_{\Delta,j}}\right)
\int_{\bfC^\sharp_{\Delta}(h^\sharp_J)} f(x) \, dx
\bvert_{h_J^\sharp=0}  \!\!\! + R^P_{w,k}(f),
\end{align}
with $h_J^\sharp = (h^\sharp_{i_1},\ldots,
h^\sharp_{i_{n_\Delta}})$ and   
\begin{equation}\label{eq:remaindersimplepoly2} R^P_{w,k}(f) :=
\sum_{\Delta \in \PB} (-1)^{m_\Delta+\dim{\Delta} } \varphi(\varepsilon,\Delta)
R^{\bfC^\sharp_{\Delta}}_{q_{J_\Delta}^\sharp, k}(f),
\end{equation}
where now the tangent cones $\bfC_\Delta$ are polarized by the
vectors $\varepsilon-\beta(\varepsilon,\Delta)$.
\end{Remark}

Let us now analyze some properties of the groups $\Gamma_\Delta$.
These generalize Claims 61, 62 and 65 in \cite{KSW} to spaces
$\Delta\in \PB$ of arbitrary dimensions. Their proofs follow easily
from the ones in \cite{KSW} but we will include them for
completeness. For that we will first introduce some necessary
notation. If $\Delta$ and $\widetilde{\Delta}$ are two elements of
$\PB$ with $\Delta \subseteq \widetilde{\Delta}$, then clearly
$\Gamma_{\widetilde{\Delta}}\subseteq \Gamma_{\Delta}$. Hence we can
define a subset $\Gamma_{\Delta}^\flat$ of $\Gamma_\Delta$ by
\begin{equation*}
\Gamma_{\Delta}^\flat:= \Gamma_{\Delta}\setminus
\bigcup_{\widetilde{\Delta} \in \PB \mid \, \Delta\subsetneq
\widetilde{\Delta}} \Gamma_{\widetilde{\Delta}}
\end{equation*}
and then
\begin{equation}\label{eq:GammaDelta}
\Gamma_\Delta = \bigsqcup_{\widetilde{\Delta}\in \PB\mid\, \Delta
\subseteq \widetilde{\Delta}} \Gamma_{\widetilde{\Delta}}^\flat 
\end{equation}
\begin{Claim}\label{claim1}
If $\gamma \in \Gamma_\Delta$ and $j \in J_\Delta$, then
$\lambda_{\gamma,j,\Delta'}$ is the same for all $\Delta'\subset
\Delta$.
\end{Claim}
\begin{Claim}\label{claim2}
If $\gamma \in \Gamma_\Delta$, $\Delta' \subset \Delta$ and $j \in
J_{\Delta'} \setminus J_\Delta$, then
$\lambda_{\gamma,j,\Delta'}=1$.
\end{Claim}
\begin{Claim}\label{claim3}
If  $\gamma \in \Gamma_\Delta^\flat$ and $j\in \Gamma_\Delta$, then
$\lambda_{\gamma,j,\Delta} \neq 1$.
\end{Claim}
\begin{proof}
Let $\gamma \in \Gamma_\Delta$ be represented by
$$
\tilde{\gamma}=\sum_{i \in J_\Delta} b_i \eta_i \in T_{J_\Delta} \cap
(\Z^d)^*
$$
for some $b_i \in \R$. Let $\Delta' \subset \Delta$. Since the 
$\alpha_{\Delta',j}$'s are dual to the  $\eta_j$'s, for $j\in
J_{\Delta'}$, we have
$$
\langle \tilde{\gamma}, \alpha_{\Delta',j} \rangle = \left\{
\begin{array}{ll} b_j, & \text{if} \quad j \in J_\Delta \\ 0, &
\text{if} \quad j  \in J_{\Delta'}\setminus J_\Delta
\end{array}\right. .
$$
Consequently,
$$
\lambda_{\gamma,j,\Delta'} = \left\{ \begin{array}{ll} e^{2 \pi i
b_j}, & \text{if} \quad j \in J_\Delta \\ 1, & \text{if} \quad j \in
J_{\Delta'}\setminus J_\Delta  \end{array} \right.
$$
is independent of $\Delta'$ and is equal to $1$ if $j \in
J_{\Delta'}\setminus J_\Delta$, and so Claim~\ref{claim1} and
\ref{claim2} follow.

Let $j \in J_\Delta$. If $\lambda_{\gamma,j,\Delta}:=e^{2\pi i
b_j}=1$, then $b_j \in \Z$ and so
\begin{equation}\label{eq:claim}
\tilde{\gamma}=\sum_{i \in J_\Delta\setminus \{j\}} b_i \eta_i
\end{equation}
also represents $\gamma$. Let $\widetilde{\Delta} \supset \Delta$ be
the element of $\PB$ such that $J_{\widetilde{\Delta}}=J_\Delta
\setminus \{j\}$. Then, by \eqref{eq:claim}, $\gamma \in
\Gamma_{\widetilde{\Delta}}$, and Claim~\ref{claim3} follows.
\end{proof}

With these properties we can further simplify formula
\eqref{eq:emlsimplepoly}. First, note that either
$h^\sharp_{\Delta,j}=h_j$, $\lambda_{\gamma,j,\Delta}^\sharp
=\lambda_{\gamma,j,\Delta}$ and $q_j^\sharp =
q_j$, or $h^\sharp_{\Delta,j}=-h_j$,
$\lambda_{\gamma,j,\Delta}^\sharp = \lambda_{\gamma,j,\Delta}^{-1}$
and $q_j^\sharp = 1 - q_j$, and so, by symmetry
property \eqref{eq:symbfN}, this gives
\begin{equation}\label{eq:symm}
{\bf N}^{k,\lambda^\sharp_{\gamma,j,\Delta}}_{q^\sharp_{\Delta,j}}
\left(\frac{\partial}{\partial h^\sharp_{\Delta,j}}\right) = {\bf
N}^{k,\lambda_{\gamma,j,\Delta}}_{q_j}\left(\frac{\partial}{\partial
h_{j}}\right).
\end{equation}
Moreover, from Claim~\ref{claim2}, we have
$\lambda_{\gamma,j,\Delta}=1$ for $j\notin J_\Delta$, implying that
$${\bf
N}^{k,\Delta_{\gamma,j,\Delta}}_{q_j}(\frac{\partial}{\partial
h_j})=1 +\text{powers of} \quad \frac{\partial}{\partial h_j}.$$
Since, still for $j\notin J_\Delta$, the cone
$\bfC^\sharp_\Delta(h_\Delta^\sharp)$ is independent of $h_j$,
\eqref{eq:emlsimplepoly} is equal to
\begin{equation}\label{eq:emlsimplepoly2}
\sum_{\Delta \in \PB} (-1)^{m_\Delta} \varphi(\varepsilon,\Delta)
\sum_{\gamma \in \Gamma_\Delta} \,\, \prod_{j=1}^N {\bf
N}^{k,\lambda_{\gamma,j,\Delta}}_{q_{j}}
\left(\frac{\partial}{\partial h_{j}}\right)
\int_{\bfC^\sharp_{\Delta}(h^\sharp_J)} f(x) \, dx \bvert_{h=0}  + R^P_{w,k}(f),
\end{equation}
where $N$ is the number of facets of $P$.
Defining
\begin{equation}\label{eq:N k,g,D}
{\bf N}_{\gamma,\Delta}^k  := \prod_{j=1}^N  {\bf
N}^{k,\lambda_{\gamma,j,\Delta}}_{q_{j}}\left(\frac{\partial}{\partial
h_{j}}\right), \quad \text{for}\quad \gamma \in \Gamma_\Delta
\end{equation}
we have, from Claim~\ref{claim1}, that
\begin{equation}\label{eq:N k,g,D'}
{\bf N}_{\gamma,\Delta}^k={\bf N}_{\gamma, \widetilde{\Delta}}^k
\quad \text{whenever} \quad  \gamma \in \Gamma_{\widetilde{\Delta}}
\quad \text{and} \quad \Delta \subset \widetilde{\Delta}.
\end{equation}
Consequently, using \eqref{eq:GammaDelta}, formula \eqref{eq:emlsimplepoly2}
can be written as
\begin{align}\label{eq:emlsimplepoly3}
 & \sum_{\Delta \in \PB}  (-1)^{m_\Delta}
\varphi(\varepsilon,\Delta) \sum_{\gamma \in \Gamma_\Delta} {\bf
N}_{\gamma,\Delta}^k \int_{\bfC^\sharp_{\Delta}
(h^\sharp_{J_\Delta})} f(x) \, dx \bvert_{h=0} +  R^P_{w,k}(f)\\ & \nonumber= 
 \sum_{\Delta \in \PB} (-1)^{m_\Delta} \varphi(\varepsilon,\Delta) \sum_{\widetilde{\Delta} \in \PB} \,\, \sum_{\gamma
\in \Gamma_{\widetilde{\Delta}}^\flat} {\bf
N}_{\gamma,\widetilde{\Delta}}^k \int_{\bfC^\sharp_{\Delta}
(h^\sharp_{J_\Delta})} f(x) \, dx \bvert_{h=0} +  R^P_{w,k}(f)
\\ \nonumber
& =  \sum_{\widetilde{\Delta} \in \PB} \,\, \sum_{\gamma
\in  \Gamma_{\widetilde{\Delta}}^\flat} {\bf
N}_{\gamma,\widetilde{\Delta}}^k \sum_{\Delta \subset
\widetilde{\Delta}} (-1)^{m_\Delta}  \varphi(\varepsilon,\Delta) 
\int_{\bfC^\sharp_{\Delta}(h^\sharp_J)} f(x) \, dx \bvert_{h=0} + R^P_{w,k}(f).
\end{align}
In the interior summation on the left we can add similar terms that correspond
to spaces $\Delta$ not included in $\widetilde{\Delta}$. indeed,
these make a zero contribution for the following reason: if $\Delta$
is not a subset of $\widetilde{\Delta}$ then there exists a $j\in
J_{\widetilde{\Delta}} \setminus J_\Delta$; then since $j \notin
J_\Delta$, the cone $\bfC^\sharp_{\Delta} (h^\sharp_{J_\Delta})$
does not depend on $h_j$; on the other hand, since $\gamma \in
\Gamma_{\widetilde{\Delta}}^\flat$ and $j\in
J_{\widetilde{\Delta}}$, we know, from Claim~\ref{claim3}, that
$\lambda_{\gamma,j,\widetilde{\Delta}}\neq 1$ and so, by
Remark~\ref{re:multipleofS}, we have that  $ {\bf
N}^{k,\lambda_{\gamma,j,\widetilde{\Delta}}}_{q_{j}}
(\frac{\partial}{\partial h_{j}})$ (one of the factors of ${\bf
N}_{\gamma,\widetilde{\Delta}}^k$) is a multiple of
$(\frac{\partial}{\partial h_{j}})$. Therefore,
\eqref{eq:emlsimplepoly3} is equal to
\begin{equation}\label{eq:emlsimplepoly4}
\begin{split}
\sum_{\widetilde{\Delta} \in \PB} \,&\, \sum_{\gamma \in
\Gamma_{\widetilde{\Delta}}^\flat} {\bf
N}_{\gamma,\widetilde{\Delta}}^k \sum_{\Delta \in \PB }
(-1)^{m_\Delta}  \varphi(\varepsilon,\Delta) 
\int_{\bfC^\sharp_{\Delta}(h^\sharp_J)} f(x) \, dx \bvert_{h=0} +
R^P_{w,k}(f) \\
& = \sum_{\widetilde{\Delta} \in \PB} \,\, \sum_{\gamma \in
\Gamma_{\widetilde{\Delta}}^\flat} {\bf
N}_{\gamma,\widetilde{\Delta}}^k  \int_{P(h)}  f(x) \, dx
\bvert_{h=0} +  R^P_{w,k}(f),
\end{split}
\end{equation} and we have our result:
\begin{Theorem}\label{thm:4}
Let $P$ be a simple polytope in $\R^d$ with $N$ facets and let
$f\in\mathcal{C}_c^{dk}(\R^d)$ be a compactly supported function on
$\R^d$ for $k\geq 1$. Choosing an $\varepsilon$ on a Paradan region
determined by $P$, we obtain
\begin{equation}\label{eq:emlfinal}
\sum_{P\cap \Z^d}{}^w f =  \sum_{\Delta \in \PB} \sum_{\gamma \in
\Gamma_{\Delta}^\flat} {\bf N}_{\gamma,\Delta }^k  \int_{P(h)}
f(x) \, dx \bvert_{h=0} +  R^P_{w,k}(f),
\end{equation}
where ${\bf N}_{\gamma,\Delta}^k$ is the  differential operator
described in \eqref{eq:N k,g,D} and \eqref{eq:N k,g,D'}, and where
the remainder is given by \eqref{eq:remaindersimplepoly}. The
operator ${\bf N}_{\gamma,\Delta}^k$ is of order $\leq k$ in each
of the variables $h_1, \ldots, h_N$ with $N$ the number of facets of
$P$. The remainder is a sum of integrals over sectors, of bounded
periodic functions times several partial derivatives of $f$ of order
no less than $k$ and no more than $k d$. Moreover, this remainder is
independent of the choice of Paradan region of $\varepsilon$, and is a distribution
supported on the polytope $P$.
\end{Theorem}

\begin{Remark}
If we instead use the polytope decompositions of Theorem~\ref{th:WBG+} the remainder in \eqref{eq:emlfinal} will be  given by \eqref{eq:remaindersimplepoly2}.
\end{Remark}

The Euler Maclaurin formula \eqref{eq:emlfinal} obtained in Theorem~\ref{thm:4} is similar to the one presented in \cite{AW}. However, in our formula, we allow the operators ${\bf
N}^{k,\lambda_{\gamma,j,\Delta}}_{q_{j}}$ that define ${\bf N}_{\gamma,\Delta}^k$ and ${\bf N}_0^{k}$  to have different weights $q_j \in \C$ while, in  \cite{AW}, the  $q_j$'s are all equal to some fixed complex number (in \cite{KSW} this fixed weight is $1/2$).  Moreover, we obtain a different expression for the remainder
$R^P_{w,k}(f)$ which is now given as a sum over the affine spaces generated by all  the faces of the polytope (not only over the vertices).
In addition,  the intermediate formulas that we obtain in \eqref{eq:emlsimplepoly3} (before adding terms with zero contribution in order to get an integral over the dilated polytope) also involve sums of integrals over the polarized tangent cones to the polytope at the different faces and not only at vertices. 

Just as the Euler Maclaurin formulas in \cite{AW}  and  \cite{KSW} our formulas  generalize to symbols\footnote{that is, smooth functions $f\in \mathcal{C}^\infty (\R^d)$ for which there is a positive integer $N$ (called the \emph{order} of the symbol) such that, for every $d$-tuple of non-negative integers $a := (a_1,\ldots ,a_d)$ there is a constant $C_a$ satisfying $\vert \partial_1^{a_1} \cdots \partial_d^{a_d} f(x) \vert \leq C_a (1+\vert x \vert)^{N-a}$.}, giving rise, in particular,  to the  following exact formula for a polynomial function $p$ in $\R^d$
\begin{equation}\label{eq:exact}
\sum_{P \cap \Z^d}{}^{w} \, p = \sum_{\Delta \in \mathcal{B} } \sum_{\gamma \in \Gamma^\flat_\Delta}  {\bf N}_{\gamma,\Delta}^k  \int_{P(h)} p(x) \, dx \bvert_{h=0}
\end{equation}
(where we choose $k \geq \text{deg} \, p + d+ 1$). From Remark~\ref{rmk:Vergne} we see that this formula is a weighted version of the exact Euler Maclaurin formula obtained in \cite{BrV}, which is obtained from \eqref{eq:exact} by making all the weights in $w$ equal to $1$.

{\bf Acknowledgments:}
We would like to thank Jonathan Weitsman for calling our attention to the
localization method of Paradan applied to the norm square of the moment map
and suggesting its use to obtain new polytope decompositions. We are also grateful to Mich\`ele Vergne for her comments on a previous version of this work.


\begin{thebibliography}{GDR}

\bibitem[A1]{A03}
J. Agapito, \emph{A weighted version of quantization commutes with reduction for a toric manifold, Integer points in polyhedra $---$ geometry, number theory,
algebra, optimization}, Contemp. Math. \textbf{374} (2005), 1--14. Proceedings of an AMS-IMS-SIAM Joint
Summer Research Conference on Integer Points in Polyhedra, July
13-17, 2003, Snowbird, Utah.

\bibitem[A2]{A05}
J. Agapito, \emph{Weighted Brianchon-Gram decomposition}. Preprint
(2004). To appear in Canad. Math. Bull.

\bibitem[At]{At}
M. Atiyah, \emph{Convexity and commuting Hamiltonians}, Bull.
London Math. Soc. \textbf{14} (1981), 1--15.

\bibitem[AB]{AB}
M. Atiyah and  R. Bott, \emph{The moment map and equivariant
cohomology}, Topology \textbf{23} (1984), no. 1, 1--28.

\bibitem[AW]{AW}
J. Agapito and J. Weitsman, \emph{The weighted Euler-Maclaurin formula
for a simple integral polytope}, Asian J. Math. \textbf{9} (2005),
no. 2, 199--212.

\bibitem[B]{B} N. Bourbaki, \emph{\'Elements de math\'ematique. XII. Premier\`ere partie: Les structures fondamentales de l'analyse. Livre IV: Functions d'une variable R\'eele}, Chapitre VI:
D\'eveloppements tayloriens g\'en\'eralis\'es; formule sommatoire
d'Euler-Maclaurin, Actualités Sci. Ind., {\bf no. 1132}, Hermann et Cie., Paris, 1951.

\bibitem[Br]{Br} C. Brianchon, \emph{Th\'eor\`{e}me nouveau sur les poly\`{e}dres}, J. \'Ecole polytechnique {\bf 15} (1837), 317-319.

\bibitem[BGV]{BGV} N. Berline, E. Getzler and M. Vergne,
\emph{Heat Kernels and Dirac Operators}, Grundlehren Math. Wiss.
\textbf{298}, Springer, Berlin (1991).

\bibitem[BV1]{BV0}
N. Berline and  M. Vergne, \emph{Z\'{e}ros d'un champ de vecteurs et classes caract\'{e}ristiques \'{e}quivariantes}, Duke Math. J. \textbf{50}
(1983), 539--549.

\bibitem[BV2]{BV1} N. Berline and M. Vergne,
\emph{Classes caract\'eristiques \'equivariantes. Formule de
localisation en cohomologie \'equivariante}, C. R. Acad. Sci. Paris
\textbf{295} (1982), 539--541.

\bibitem[BV3]{BV2} N. Berline and M. Vergne, \emph{Local Euler-Maclaurin
formula for polytopes}, Preprint (2005): arXiv:math.CO/0507256v1.

\bibitem[BrV1]{BrV} M. Brion and M. Vergne, \emph{Lattice points in simple polytopes}, J. Amer. Math. Soc. {\bf 10}, no. 2, (1997), 371-392.

\bibitem[BrV2]{BrV2} M. Brion and M. Vergne, \emph{Residue formulae, vector partition functions and lattice points in rational polytopes}, J. Amer. Math. Soc. {\bf 10}, no. 4, (1997), 797-833.

\bibitem[C]{C}
J. Cassels, \emph{An Introduction to the Geometry of Numbers}, Classics in Mathematics, Springer (1997).


\bibitem[CS1]{CS1} S. Cappell and J. Shaneson, \emph{Genera of algebraic varieties and counting lattice points},  Bull. Amer. Math. Soc. (N.S.) {\bf 30} (1994), 62-69.

\bibitem[CS2]{CS2} S. Cappell and J. Shaneson, \emph{Euler Mac-Laurin expansions for lattices above dimension one}, C. R. Acad. Sci. Paris S\'{e}r A {\bf 321} (1995), 885-890.


\bibitem[De]{De} T. Delzant, \emph{Hamiltoniens p\'eriodiques et image
convexe de l'application moment},   Bull. Soc. Math. France \textbf{116} (1988), 315--339.

\bibitem[DV]{DV} M. Duflo and M. Vergne, \emph{Cohomologie \'equivariante et descente}, Ast\'erisque
\textbf{215} (1993), 5--108.

\bibitem[G]{G} J. Gram, \emph{Om rumvinklerne i et polyeder}, Tidsskrift for Math. (Copenhagen) {\bf 4} no 3. (1874), 161-163. 

\bibitem[Gr]{Gr} B. Gr\"{u}nbaum, \emph{Complex polytopes}, Second Edition Prepared by Volker Kaibel, Victor Klee, and G\"{u}nter Ziegler, Graduate Texts in Mathematics {\bf 221}, Springer-Verlag, New York (2003). 


\bibitem[Gu1]{Gu1}
V. Guillemin, \emph{Moment maps and Combinatorial Invariants of
Hamiltonian $T^n$-spaces}, Progress in Mathematics \textbf{122},
Birkh\"auser (1994).

\bibitem[Gu2]{Gu2}
V. Guillemin, \emph{Riemann-Roch for toric orbifolds}, J. Differ. Geom. \textbf{45} (1997), 53--73.


\bibitem[GGK]{GGK} V. Guillemin, V. Ginzburg and Y. Karshon, \emph{Moment maps, cobordisms, and Hamiltonian group actions},
Mathematical Surveys and Monographs {\bf 98}, American Mathematical Society, Providence, RI (2202).

\bibitem[GS]{GS}
V. Guillemin and S. Sternberg, \emph{Convexity properties of the moment map}, Invent. Math. \textbf{67} (1982), 491--513.

\bibitem[H]{H} L. Hormander, \emph{Lectures on non-linear hyperbolic differential equations}, Math\'{e}matiques \& Applications (Berlin), {\bf 26}, Springer-Verlag, Berlin, 1997. 

\bibitem[K]{K} F. Kirwan, \emph{Cohomology of quotients in
symplectic and algebraic geometry}, Princeton Univ. Press, Princeton
(1984).

\bibitem[Kn]{Kn} K. Knopp, \emph{Theory and Application of Infinite Series}, Dover Publications, New York (1990). First published in German in 1921 and in English in 1928.

\bibitem[KP1]{KP1} A. Khovanskii and A. Pukhlikov, \emph{Finitely additive measures of virtual polytopes}, (Russian)  Algebra i Analiz {\bf 4} (1992), no. 2, 161-185; translation in St. Petersburg Math. J. {\bf 4} (1993), no. 2, 337-356.

\bibitem[KP2]{KP2} A. Khovanskii and A. Pukhlikov, \emph{The Riemann-Roch theorem for integrals and sums of quasipolynomials on virtual polytopes}, Algebra i Analiz {\bf 4} (1992), no. 4,  188-216; translation in St. Petersburg Math. J. {\bf 4} (1993), no. 4, 789-812.

\bibitem[KSW]{KSW} Y. Karshon, S. Sternberg and J. Weitsman, \emph{Euler-
Maclaurin with remainder for a simple integral polytope}, Preprint
(2003): arXiv:math.CO/0307125v1.

\bibitem[KV]{KV} S. Kumar and M. Vergne, \emph{Equivariant cohomology
with generalized coefficients}, Ast\'erisque \textbf{215} (1993),
109--204.

\bibitem[L]{L} J. Lawrence, \emph{Polytope volume computation}, Math. Comp {\bf 57} (1991), no. 195, 259-271.

\bibitem[LT]{LT} E. Lerman and S. Tolman, \emph{Hamiltonian torus actions on
symplectic orbifolds and toric varieties}, Trans. Amer. Math. Soc.
\textbf{349} (1997), no. 10, 4201--4230.

\bibitem[M]{M}
R. Morelli, \emph{Pick's theorem and the Todd class of a toric variety}, Adv. Math.
\textbf{100} (1993), no. 2, 183--231.

\bibitem[P]{P}
P-E. Paradan, \emph{Formules de localisation en cohomologie
equivariante}, Compositio Math \textbf{117} (1999), 243--293.

\bibitem[Po]{Po}
J. Pommersheim, \emph{Toric varieties, lattice points and Dedekind sums}, Math. Ann. \textbf{295} (1993), 1--24.

\bibitem[PT]{PT}
J. Pommersheim and H. Thomas, \emph{Cycles representing the the Todd class of a toric variety},  J. Amer. Math. Soc.  {\bf 17}  (2004),  no. 4, 983--994 (electronic).

\bibitem[S]{S} J. Shaneson, \emph{Characteristic Classes, Lattice Points, and Euler-Maclaurin Formulae}, Proceedings of the International Congress of Mathematicians, Zurich, 1994. Basel: Birkauser Verlag, 1995.

\bibitem[So]{So} D. Sommerville, \emph{the relations connecting the angle-sums and volume of a polytope in space of $n$ dimensions}, Proceedings Royal Soc. London (A) {\bf 115} (1927), 103--119.

\bibitem[SV]{SV} A. Szenes and M. Vergne, \emph{Residue formulae for vector partitions and Euler-MacLaurin sums}, Adv. in Appl. Math. {\bf 30} (2003), no. 1-2, 295-342.

\bibitem[V]{V} A. Varchenko, \emph{Combinatorics and topology of the arrangement of affine hyperplanes in the real space} (Russian) Funktsional. Anal. i Prilozhen. {\bf 21}, no. 1, 11--22. English translation: Functional Anal. Appl. {\bf 21} (1987), no. 1, 9--19.
\end{thebibliography}
\end{document}